\documentclass[openright,11pt,two0side,a4paper,leqno]{article}
\usepackage[latin1]{inputenc} 
\usepackage{amsmath,amssymb} 
\usepackage{theorem}

\newlength{\larg}
\newcommand{\bd}[1]{\boldsymbol #1}
\newcommand{\ds}{\displaystyle}

\newcommand{\pv}{\noindent \emph{Proof. }}
\newcommand{\cqfd}{\phantom{m} \hfill $\Box$}
\newcommand{\fini}{\phantom{m} \hfill $\diamond$}
\newenvironment{proof}{\pv}{\cqfd \\ }

\newcommand{\ie}{\emph{i.e. }}
\newcommand{\eg}{\emph{e.g. }}
\newcommand{\resp}{\emph{resp. }}

\newcommand{\barre}{\overline{ ^{\ \: }} \ }

\newcommand{\set}[1]{ \{ #1 \}}

\newcommand{\mN}{{\mathbb N}}
\newcommand{\mQ}{{\mathbb Q}}
\newcommand{\mR}{{\mathbb R}}
\newcommand{\mZ}{{\mathbb Z}}

\newcommand{\interv}[2]{[\! [ #1  ;  #2 ] \!]}

\newcommand{\wt}{\mathrm{wt}}
\newcommand{\wtpt}{\dot{\mathrm{wt}}}

\newcommand{\eilow}{\tilde{e_i}^{\scriptsize \mbox{\rm low}}}

\newcommand{\filow}{\tilde{f_i}^{\scriptsize \mbox{\rm low}}}

\newcommand{\gsl}[1]{{\mathfrak{sl}}_{#1}}
\newcommand{\hsl}[1]{{\widehat{\mathfrak{sl}}}_{#1}}
\newcommand{\Uq}{U_q(\hsl{n})}
\newcommand{\Up}{U_p(\hsl{l})}
\newcommand{\Uqprime}{U_q'(\hsl{n})}
\newcommand{\Upprime}{U_p'(\hsl{l})}
\newcommand{\Mq}[1][\bd{s}_l]{\boldsymbol{\operatorname{M}}_q[#1]}
\newcommand{\Fq}[1][\bd{s}_l]{\boldsymbol{\operatorname{F}}_q[#1]}
\newcommand{\Fp}[1][\bd{s}_n]{\boldsymbol{\operatorname{F}}_p[#1]^{\bullet}}

\newtheorem{thm}{Theorem}[section]
\newtheorem{lemma}[thm]{Lemma}
\newtheorem{prop}[thm]{Proposition}
\newtheorem{cor}[thm]{Corollary}

{\theorembodyfont{\rmfamily} \newtheorem{remark}[thm]{Remark}}   %supprime l'italique dans cet environnement \newtheorem
{\theorembodyfont{\rmfamily} \newtheorem{example}[thm]{Example}} % (nécessite le package 'theorem')
{\theorembodyfont{\rmfamily} \newtheorem{definition}[thm]{Definition}}
{\theorembodyfont{\rmfamily} \newtheorem{notation}[thm]{Notation}}

\title{An algorithm for computing the canonical bases of higher-level $q$-deformed Fock spaces}
\author{Xavier YVONNE}
\date{}

\begin{document}

\maketitle

\begin{abstract} We derive a straightening-free algorithm that computes the canonical bases of any higher-level $q$-deformed Fock space.
\end{abstract}


\section{Introduction} \label{section_intro}

The higher-level $q$-deformed Fock spaces form an important family of integrable representations of the quantum group $\Uq$. 
The Fock representation $\Fq$ depends on a parameter $\bd{s}_l=(s_1,\ldots,s_l) \in \mZ^l$ called multi-charge. It was introduced in \cite{JMMO} 
in order to compute the crystal graph of the integrable irreducible representation of $\Uq$ with highest weight 
$\Lambda_{s_1}+\cdots+\Lambda_{s_l}$.

The canonical basis of a Fock space is a lower global crystal basis in the sense of \cite{Kas1} that enjoys nice properties. It was constructed for 
$l=1$ in
\cite{LT1,LT2} and for $l \geq 1$ by Uglov \cite{U}. This construction relies on 
\begin{enumerate}
\item an embedding of $\Uq$-modules: $\Fq \hookrightarrow \Lambda^s$ which is a bijection if (and only if) $l=1$. Here $\Lambda^s$ is the space
of semi-infinite $q$-deformed wedge products of charge $s=s_1+\cdots+s_l$ (see \cite{U}).
\item the definition of a $\barre$-involution of $\Lambda^s$ compatible with the embedding above. The canonical basis is then the unique 
$\barre$-invariant basis of $\Fq$ satisfying a certain congruence property modulo the $\mZ[q]$-lattice spanned by $l$-multi-partitions.
\end{enumerate} 
In \cite{U}, Uglov provides an algorithm for computing the canonical basis based on the straightening of non-ordered $q$-wedge products. 
Unfortunately, in particular when $s_1 \gg \cdots \gg s_l$, the number of $q$-wedge products (of $2$ factors) to be straightened becomes too large for this
algorithm being used for practical computations.

The goal of this paper is to derive a faster algorithm, which does not require straightening $q$-wedge products. First, we compute a 
$\barre$-invariant basis of $\Lambda^s$. The vectors of this basis are obtained by letting act the generators 
of the quantum groups $\Uq$ and $\Up$ (with $p=-q^{-1}$) and a Heisenberg algebra $\mathcal{H}$
to the vacuum vector (see Thm. \ref{thm_base_Lambda_s_barre_invariante}). Then, we compute the intersection of this basis with a given weight subspace $V$ of $\Fq$ 
(see Prop. \ref{prop_base_sep_de_Fqsl_barre_invariante}). We thus get the transition matrix $T(q)$ between the standard basis and a $\barre$-invariant 
basis of $V$. We can then compute the matrix $A(q)$ of the involution of $V$ with respect to the standard basis by the formula
$A(q)=T(q)\big(T(q^{-1})\big)^{-1}$. Note that we are forced here to invert a large matrix with Laurent polynomial entries. This is the time-consuming step of our algorithm. Note also
that although $A(q)$ is unitriangular by \cite{U}, this matrix $T(q^{-1})$ is not. Once $A(q)$ is known, the computation of the canonical basis is pure routine (see \eg \cite[Thm. 7.1]{L}). 

Our algorithm is a generalization of the algorithm of Leclerc and Thibon (see \cite{L}) for $l=1$. However, two difficulties arise when
$l>1$. First, we have to take into account the action of $\Up$, which is trivial if $l=1$. Secondly, the action of the Heisenberg algebra 
$\mathcal{H}$ is much more complicated for $l>1$ and it particular it requires some new work to calculate its
action on the vacuum vectors without using the straightening relations (see Section \ref{section_action_H}).  

Using our algorithm, we were able to compute the canonical bases of large weight subspaces of different Fock spaces (see examples in 
Section \ref{section_tables}). These calculations helped us to conjecture and prove a theorem giving a combinatorial expression of the derivative at 
$q=1$ of the $\barre$-involution (see \cite[Thm. 2.11]{Y1}). The latter result supports in turn a new conjecture for computing the $q$-decomposition
matrices of Dipper-James-Mathas' cyclotomic $v$-Schur algebras with parameters specialized at powers of a complex $n$-th root of unity 
(\cite{DJM,Y1}). In this context, our result is an analogue of the Jantzen sum formula for cyclotomic $v$-Schur algebras \cite{JM}.\\       

\noindent \emph{Acknowledgments.} I am grateful to Bernard Leclerc for the illuminating discussions I 
had with him. I thank him for his reading of the preliminary versions of this paper. I also thank the referee
for his/her very quick replies.\\

\noindent \emph{Notation.} Let $\mN$ (\resp $\mN^{*}$) denote the set of nonnegative (\resp positive) integers, and for $a,\,b \in \mR$ denote by 
$\interv{a}{b}$ the discrete interval $[a  ;  b] \cap \mZ$. For $X \subset \mR$, $a \in \mR$, $N \in \mN^*$, put 
\begin{equation} \label{eq_Z_l_s}
X^N(a):=\set{(x_1,\ldots,x_N) \in X^N \mid x_1+\cdots+x_N=a}.
\end{equation}
Throughout this article, we fix 3 integers $n$, $l \geq 1$ and $s \in \mZ$. Let $\Pi$
denote the set of all integer partitions. If $(W,S)$ is a Coxeter system, denote by $\ell : W \rightarrow \mN$ the
length function on $W$.

\section{Higher-level $q$-deformed Fock spaces}
In this section, we recall briefly the definition of the higher-level Fock spaces and their canonical bases.
These objects were introduced by \cite{U}, to which we refer the reader for more details. We follow here the
notation from \cite{Y2}. 

\subsection{The quantum algebras $\Uq$ and $\Up$}

In this section, we assume that $n \geq 2$ and $l \geq 2$. Let $\hsl{n}$ be the Kac-Moody algebra of type $A_{n-1}^{(1)}$ defined over the field 
$\mQ$ \cite{Kac}. Let $\Lambda_0,\ldots,\Lambda_{n-1}$ be the fundamental weights, $\alpha_0,\ldots,\alpha_{n-1}$ be the simple roots and 
$\delta:=\alpha_0+\cdots+\alpha_{n-1}$ be the null root. It will be convenient to extend the index set of the fundamental weights by setting 
$\Lambda_i:=\Lambda_{i \bmod n}$ for all $i \in \mZ$. The space 
$\bigoplus_{i=0}^{n-1}{\mQ\,\Lambda_i} \; \oplus \; \mQ\,\delta = \bigoplus_{i=0}^{n-1}{\mQ\,\alpha_i} \; \oplus \; \mQ\,\Lambda_0$ 
is equipped with a non-degenerate bilinear symmetric form $(.\,,\,.)$ defined by
\begin{equation}
(\alpha_i,\alpha_j)=a_{i,j}, \quad (\Lambda_0,\alpha_i) = \delta_{i,0}, \quad 
(\Lambda_0,\Lambda_0)=0 \qquad (0 \leq i,j \leq n-1),
\end{equation}
where $(a_{i,j})_{0 \leq i,j \leq n-1}$ is the Cartan matrix of $\hsl{n}$. Let $\Uq$ be the $q$-deformed universal enveloping algebra of $\hsl{n}$ 
(see \eg \cite{KMS,U}). This is an algebra over $\mQ(q)$ with generators $e_i$, $f_i$, $t_i^{\pm 1}$ $(0 \leq i \leq n-1)$ and $\partial$; the relations 
will be omitted. The subalgebra of $\Uq$ generated by $e_i$, $f_i$, $t_i^{\pm 1}$ $(0 \leq i \leq n-1)$ will be denoted by $\Uqprime$. If $M$ is a 
$\Uq$-module, denote by $\mathcal{P}(M)$ the set of weights of $M$ and let $M \langle w \rangle$ denote the subspace of $M$ of weight $w$. Let
\begin{equation}  
\wt(x):=w
\end{equation}
denote the weight of $x \in M \langle w \rangle$. The Weyl group of $\hsl{n}$ (or $\Uq$) is
\begin{equation}  
W_n=\langle \sigma_0,\ldots,\sigma_{n-1} \rangle \cong \widetilde{\mathfrak{S}}_n,
\end{equation}
where $\sigma_i$ $(0 \leq i \leq n-1)$ is the orthogonal reflection that fixes pointwise the hyperplane orthogonal to
$\alpha_i$. \\

We also introduce the algebra $\Up$ with
\begin{equation}
p:=-q^{-1}.
\end{equation}
In order to distinguish the elements related to $\Uq$ from those related to $\Up$, we put dots over the latter. For example, $\dot{e}_i$, 
$\dot{f}_i$, $\dot{t}_i^{\pm 1}$ $(0 \leq i \leq l-1)$ and $\dot{\partial}$ are the generators, $\dot{\alpha}_i$ $(0 \leq i \leq l-1)$ 
are the simple roots, $\dot{W}_l = \langle \dot{\sigma}_0,\ldots,\dot{\sigma}_{l-1} \rangle$ 
is the Weyl group of $\Up$ and so on.

\subsection{The space $\Lambda^s$}

Following \cite{U}, let 
\begin{equation} \label{notation_Lambdas}
\Lambda^s=\Lambda^s[n,l]
\end{equation}
be the space of (semi-infinite) $q$-wedge products of charge $s$ 
(this space is denoted by $\Lambda^{s+\frac{\infty}{2}}$ in \cite{U}). This vector space has a natural basis
formed by the ordered $q$-wedge products; this basis is called \emph{standard}. A non-ordered $q$-wedge product
can be expressed as a linear combination of ordered $q$-wedge products by using the straightening relations given
in \cite[Prop. 3.16]{U} (we do not need these relations in this article). Recall the notation $\mZ^l(s)$ and 
$\mZ^n(s)$ from (\ref{eq_Z_l_s}). Following [U, \S 4.1] or  [Y2, \S 2.2.1], 
%!!!!!!!!!!!!!!!!!!!!!!!!!!!!!!!!!!!!!!!!!!!!!!!! couldn't use \S inside the 'cite' macro
 we shall use in the sequel the following
indexations of the standard basis:
\begin{equation} 
\set{|\lambda,s \rangle \mid \lambda \in \Pi}
= \set{|\bd{\lambda}_l,\bd{s}_l \rangle \mid \bd{\lambda}_l \in \Pi^l,\,\bd{s}_l \in \mZ^l(s)}
= \set{|\bd{\lambda}_n,\bd{s}_n \rangle^{\bullet} \mid \bd{\lambda}_n \in \Pi^n,\,\bd{s}_n \in \mZ^n(s)}.
\end{equation}
Following \cite{JMMO,FLOTW,U}, the vector space $\Lambda^s$ can be made into an integrable representation of level $l$ of the quantum 
algebra $\Uq$. The action of $\Uq$ can be described in a nice way in terms of addable/removable $i$-nodes of $l$-multi-partitions;
see \cite[Eq. (33-34)]{U}. Note that these formulas do no involve the straightening of $q$-wedge products; they are therefore handy to 
use for computations. In a completely similar way, $\Lambda^s$ can be made into an integrable representation of level $n$ of the quantum algebra 
$\Up$. This action can be described using the indexation by $n$-multi-partitions; see \cite[Eq. (35-36)]{U}. The vectors of the standard
basis of $\Lambda^s$ are weight vectors for the actions of $\Uq$ and $\Up$. In order to recall the expression of the weights, let us introduce the 
following notation. For $\bd{\lambda}_l \in \Pi^l$, $\bd{s}_l \in \mZ^l(s)$, $0 \leq i \leq n-1$, denote by $N_i(\bd{\lambda}_l;\bd{s}_l;n)$ the 
number of nodes with residue $i$ (modulo $n$) that are contained in the Young diagram of $\bd{\lambda}_l$ (note that the definition of residues 
involves the multi-charge $\bd{s}_l$). For $\bd{s}_l = (s_1,\ldots,s_l) \in \mZ^l$, define
\begin{equation}
\Delta(\bd{s}_l,n):= \ds \frac{1}{2} \sum_{b=1}^{l}{\Bigl( \frac{s_b^2}{n}-s_b \Bigr) - 
\Bigl(\frac{(s_b \bmod n)^2}{n}-(s_b \bmod n)\Bigr)},
\end{equation}
where for $1 \leq b \leq n$, $s_b \bmod n$ denotes the integer in $\interv{0}{n-1}$ that is congruent to $s_b$ modulo~$n$. 

\begin{prop}[\cite{U}] \label{prop_formulas_weights} With the notation above, we have
\begin{eqnarray}
\hspace{1cm} \wt(|\bd{\lambda}_l, \bd{s}_l \rangle) &=&
 - \Delta(\bd{s}_l,n) \delta + \Lambda_{s_1}+ \cdots + \Lambda_{s_l} - \ds \sum_{i=0}^{n-1} N_i(\bd{\lambda}_l;\bd{s}_l;n) \,\alpha_i, 
\label{poidsl} \\
\hspace{1cm} \dot{\wt}(|\bd{\lambda}_l, \bd{s}_l \rangle) &=& 
 - \bigl(\Delta(\bd{s}_l,n)+ N_0(\bd{\lambda}_l;\bd{s}_l;n) \bigr) \dot{\delta} 
 + (n-s_1+s_l) \dot{\Lambda}_0+ \ds \sum_{i=1}^{l-1}{(s_i-s_{i+1})  \,\dot{\Lambda}_i},
\label{poidsptl} \\
\hspace{1cm} \dot{\wt}(|\bd{\lambda}_n, \bd{s}_n \rangle^{\bullet}) &=& 
- \Delta(\bd{s}_n,l) \dot{\delta} + \dot{\Lambda}_{s_1}+ \cdots + \dot{\Lambda}_{s_n} - \ds \sum_{i=0}^{l-1} N_i(\bd{\lambda}_n;\bd{s}_n;l) 
\,\dot{\alpha}_i, 
\label{poidsptn} \\
\hspace{1cm} \wt (|\bd{\lambda}_n, \bd{s}_n \rangle^{\bullet}) &=& 
 - \bigl(\Delta(\bd{s}_n,l)+ N_0(\bd{\lambda}_n;\bd{s}_n;l) \bigr) \delta 
  + (l-s_1+s_n) \Lambda_0 + \ds \sum_{i=1}^{n-1}{(s_i-s_{i+1}) \,\Lambda_i}. 
\label{poidsn}
\end{eqnarray} \cqfd
\end{prop}

For $m \in \mZ^*$, Uglov defined an endomorphism $B_m$ of $\Lambda^s$ (see \cite[Eq. (25) \& Sect. 4.3]{U} or
\cite[Def. 2.2]{Y2}). %réf pour la déf dans Y2!!!!!!!!!!!!!!!!!!!!!!!!!
His definition is obtained by taking the limit $r \rightarrow \infty$ in the action of the center of the Hecke 
algebra of $\widehat{\mathfrak{S}}_r$ on $q$-wedge products of $r$ factors. However, by \cite{U} the operators 
$B_m$ do not commute, but they span a Heisenberg algebra
\begin{equation}
\mathcal{H}=\langle B_m \mid m \in \mZ^* \rangle.
\end{equation}

Note that the $q$-wedge products involved in the definition of $B_m$ are in general not ordered. Therefore, the 
computation of the action of $\mathcal{H}$ often requires straightening many $q$-wedge products. \\

We now recall some results concerning the actions of $\Uq$, $\Up$ and $\mathcal{H}$ on $\Lambda^s$. 

\begin{prop}[\cite{U}] The actions of $\Uqprime$, $\Upprime$ and $\mathcal{H}$ on $\Lambda^s$ pairwise commute.\hspace{-3mm}\cqfd
\end{prop} 

For $L$, $N \in \mN^*$, introduce the finite set
\begin{equation} \label{definition_A_l_n_s}
A_{L,N}(s):= \set{(r_1,\ldots,r_L) \in \mZ^L(s) \mid r_1 \geq \cdots \geq r_L, \ r_1-r_L \leq N}.
\end{equation}

\begin{thm}[\cite{U}, Thm. 4.8] \label{thm_dec_Lambda_s} We have
$$
\Lambda^s = \ds \bigoplus_{\bd{r}_l \in A_{l,n}(s)}{\Uqprime \otimes \mathcal{H} \otimes \Upprime.  
|\bd{\emptyset}_l,\bd{r}_l} \rangle
= \ds \bigoplus_{\bd{r}_n \in A_{n,l}(s)}{\Uqprime \otimes \mathcal{H} \otimes \Upprime.  
|\bd{\emptyset}_n,\bd{r}_n \rangle^{\bullet}}.$$ \cqfd
\end{thm}

\subsection{Higher-level $q$-deformed Fock spaces}

Recall that $p=-q^{-1}$. For $\bd{s}_l \in \mZ^l(s)$, $\bd{s}_n \in \mZ^n(s)$, let
\begin{equation}
\Fq=\ds \bigoplus_{\bd{\lambda}_l \in \Pi^l}{\mQ(q) |\bd{\lambda}_l, \bd{s}_l \rangle}, \qquad
\Fp= \ds \bigoplus_{\bd{\lambda}_n \in \Pi^n}{\mQ(q) |\bd{\lambda}_n, \bd{s}_n \rangle}^{\bullet}
\end{equation}
denote the higher-level ($q$-deformed) Fock spaces \cite{U}. By \cite{U}, the subspace $\Fq \subset \Lambda^s$ 
is stable under the actions of $\Uq$ and $\mathcal{H}$ (but not under the action of $\Up$) 
and the subspace $\Fp \subset \Lambda^s$ is stable under the actions of $\Up$ and $\mathcal{H}$. \\

For $N$, $L \in \mN^*$, define a bijective map
\begin{equation}
\theta_{L,N} \, : \mQ^{L}(s) \rightarrow \mQ^{L}(N), \quad (s_1,\ldots,s_L) \mapsto (N-s_1+s_L,s_1-s_2,\ldots, s_{L-1}-s_L).
\end{equation}
The next result shows that the Fock spaces are sums of certain weight subspaces of $\Lambda^s$.
The proof follows easily from Proposition \ref{prop_formulas_weights}.

\begin{prop}[\cite{U}] \ \\[-5mm] \label{fock_spaces_weight_subspaces}
\begin{itemize}
\item[\rm (i)] Let $\bd{s}_n \in \mZ^n(s)$. Let $(a_0,\ldots,a_{n-1}):=\theta_{n,l}(\bd{s}_n)$ and $w := \sum_{i=0}^{n-1}{a_i \Lambda}_i$.
Then we have $\Fp= \ds \bigoplus_{d \in \mZ}{\Lambda^s \langle w+d \delta} \rangle.$
\item[\rm (ii)] Let $\bd{s}_l \in \mZ^l(s)$. Let $(a_0,\ldots,a_{l-1}):=\theta_{l,n}(\bd{s}_l)$ and $\dot{w}:=\sum_{i=0}^{l-1}{a_i \dot{\Lambda}_i}$.
Then we have $\Fq= \ds \bigoplus_{d \in \mZ}{\Lambda^s \langle {\dot{w}+d \dot{\delta}}} \rangle.$ \cqfd
\end{itemize}
\end{prop}

We now compare some weight subspaces of the Fock spaces. The proof again follows from Proposition \ref{prop_formulas_weights}. 

\begin{prop} \label{prop_bij_pds_charges}
Let $\bd{s}_l \in \mZ^l(s)$ and $w$ be a weight of $\Fq$. Then there exists a unique pair $(\bd{s}_n,\dot{w})$ such that
$\bd{s}_n \in \mZ^n(s)$, $\dot{w}$ is a weight of $\Fp$ and $\Fq \langle w \rangle = \Fp \langle \dot{w} \rangle$. More precisely, write
$w = d \delta + \sum_{i=0}^{n-1}{a_i \Lambda_i}$ $($with $a_0,\ldots,a_{n-1},d \in \mZ)$, 
$\bd{s}_l=(s_1,\ldots,s_l)$ and put $s_0:=n+s_l$. Then we have 
$\bd{s}_n=\theta_{n,l}^{-1}(a_0,\ldots,a_{n-1})$ and $\dot{w}= d \dot{\delta} + \sum_{i=0}^{l-1}{(s_i-s_{i+1}) \dot{\Lambda}_i}.$ \cqfd
\end{prop}

\begin{example} Take $n=3$, $l=2$, $\bd{s}_l=(1,0)$ and $w = -2 \Lambda_0 + \Lambda_1 + 3 \Lambda_2 - 2 \delta$. Then by (\ref{poidsl}),
we have $\wt \big(\big|\bigl((1,1),(1)\bigr),\bd{s}_l \big\rangle \big)=w$, so $w$ is a weight of $\Fq$.
By Proposition \ref{prop_bij_pds_charges}, we have $\Fq \langle w \rangle = \Fp \langle \dot{w} \rangle$ with $\bd{s}_n=(2,1,-2)$ and 
$\dot{w} = 2 \dot{\Lambda}_0+\dot{\Lambda}_1-2 \dot{\delta}$. Moreover, using
(\ref{poidsl}) and (\ref{poidsptn}), we see that for all 
$|\bd{\lambda}_l,\bd{s}_l \rangle = |\bd{\lambda}_n,\bd{s}_n \rangle^{\bullet} \in \Fq \langle w \rangle = \Fp \langle \dot{w} \rangle$,
we have $N_0(\bd{\lambda}_l;\bd{s}_l;n)=2$, $N_1(\bd{\lambda}_l;\bd{s}_l;n)=1$, $N_2(\bd{\lambda}_l;\bd{s}_l;n)=0$ and
$N_0(\bd{\lambda}_n;\bd{s}_n;l)=N_1(\bd{\lambda}_n;\bd{s}_n;l)=0$ (this shows \emph{a posteriori} that 
$\operatorname{dim}(\Fq \langle w \rangle)=1$). \fini
\end{example}

\subsection{Action of the Weyl groups $W_n$ and $\dot{W}_l$}

The Weyl group $W_n$ acts on the weight lattice $\bigoplus_{i=0}^{n-1}{\mZ \Lambda_i} \oplus \mZ \delta$ by
\begin{equation} \label{eq_action_Weyl_group_fundamental_weights}
\sigma_i.\delta = \delta \quad \mbox{and} \quad
 \sigma_i.\Lambda_j= \left\{
 \begin{array}{lc}
 \Lambda_j & \mbox{ if } j \neq i, \\
 \Lambda_{i-1}+ \Lambda_{i+1} - \Lambda_i -\delta_{i,0} \, \delta & \mbox{ if } j = i \\
\end{array} \right. \qquad (0 \leq i,j \leq n-1).
\end{equation}
 
Moreover, it is easy to see that $W_n$ acts faithfully on $\mZ^n(s)$ by
\begin{equation}
\left\{ \begin{array}{rcl}
\sigma_0.(s_1,\ldots,s_n) &=& (s_n+l,s_2,\ldots,s_{n-1},s_1-l), \\
\sigma_i.(s_1,\ldots,s_n) &=& (s_1,\ldots,s_{i+1},s_i,\ldots,s_n) \qquad (1 \leq i \leq n-1),\\
\end{array} \right.
\end{equation}
and the set $A_{n,l}(s)$ defined by (\ref{definition_A_l_n_s}) is a fundamental domain for this action. In a similar way, one can define two actions 
of the Weyl group $\dot{W}_l$ of $\Up$, one on the weight lattice $\bigoplus_{i=0}^{l-1}{\mZ \dot{\Lambda}_i} \oplus \mZ \dot{\delta}$ and one on 
$\mZ^l(s)$.

\subsection{The lower crystal basis $(\mathcal{L}[\bd{s}_l],\mathcal{B}[\bd{s}_l])$ of $\Fq$ at $q=0$}

Let $\bd{s}_l \in \mZ^l(s)$. Following \cite{Kas1}, let $\mathbb{A} \subset \mQ(q)$ be the ring of rational functions which are regular at $q=0$, 
$\mathcal{L}[\bd{s}_l]:=\bigoplus_{\bd{\lambda}_l \in \Pi^l}{\mathbb{A} \, |\bd{\lambda}_l,\bd{s}_l \rangle}$ and
for $0 \leq i \leq n-1$, let $\tilde{e}_i=\eilow$ and $\tilde{f}_i=\filow$ denote Kashiwara's operators acting on 
$\mathcal{L}[\bd{s}_l]$. Put
\begin{equation}
\mathcal{B}[\bd{s}_l]:=\set{|\bd{\lambda}_l,\bd{s}_l \rangle \bmod q \mathcal{L}[\bd{s}_l] \mid \bd{\lambda}_l \in \Pi^l}.
\end{equation} 
In the sequel, if $\bd{s}_l \in \mZ^l(s)$ is fixed, we shall write more briefly $\bd{\lambda}_l$ for the
element in $\mathcal{B}[\bd{s}_l]$ indexed by the corresponding multi-partition. By \cite{JMMO},\cite{FLOTW},\cite{U}, the pair 
$(\mathcal{L}[\bd{s}_l],\mathcal{B}[\bd{s}_l])$ is a lower crystal basis of $\Fq$ at $q=0$ in the sense of \cite{Kas1}, and the
crystal graph $\mathcal{B}[\bd{s}_l]$ 
contains the arrow $\bd{\lambda}_l \stackrel{i}{\longrightarrow} \bd{\mu}_l$ if and only if the multi-partition $\bd{\mu}_l$ is obtained 
from $\bd{\lambda}_l$ by adding a good $i$-node in the sense of \cite[Thm. 2.4]{U}. \\

We now recall the definition of the involution $\sigma_i$ of $\mathcal{B}[\bd{s}_l]$ (we sometimes view $\sigma_i$ as a bijection of
$\Pi^l$). First, let us introduce a piece of notation that will be used in the sequel.

\begin{notation} For $\bd{s}_l \in \mZ^l(s)$ and $w \in \mathcal{P}(\Fq)$, put
\begin{equation}
\Pi^l(\bd{s}_l;w):= \set{ \bd{\lambda}_l \in \Pi^l \ \big| \ |\bd{\lambda}_l,\bd{s}_l \rangle \in \Lambda^s \langle w \rangle},
\end{equation}
and define similarly $\Pi^n(\bd{s}_n;\dot{w})$ for $\bd{s}_n \in \mZ^n(s)$ and $\dot{w} \in \dot{\mathcal{P}}(\Fp)$. \fini
\end{notation}

\begin{definition} \label{def_sigma_i} 
Fix $\bd{s}_l \in \mZ^l(s)$. Let $\bd{\lambda}_l \in \mathcal{B}[\bd{s}_l] \cong \Pi^l$ and $i \in \interv{0}{n-1}$. Let $\mathcal{C}$ be the 
$i$-chain in $\mathcal{B}[\bd{s}_l]$ containing $\bd{\lambda}_l$. Let $\sigma_i(\bd{\lambda}_l) \in \mathcal{B}[\bd{s}_l] \cong \Pi^l$ be the unique 
element in $\mathcal{C}$ such that $\wt(\sigma_i(\bd{\lambda}_l))=\sigma_i.(\wt(\bd{\lambda}_l))$. In other words, $\sigma_i(\bd{\lambda}_l)$ is 
obtained from $\bd{\lambda}_l$ via a central symmetry in the middle of $\mathcal{C}$. This defines an involution $\sigma_i$ of 
$\mathcal{B}[\bd{s}_l]$. This map induces, for $w \in \mathcal{P}(\Fq)$, a bijection 
\begin{equation} \label{def_sigma_i_2}
\sigma_i : \Pi^l(\bd{s}_l;w) \stackrel{\sim}{\longrightarrow} \Pi^l(\bd{s}_l;\sigma_i.w).
\end{equation} \fini 
\end{definition}

By \cite{Kas2}, the definition of $\sigma_0,\ldots,\sigma_{n-1}$ as bijections of 
$\mathcal{B}[\bd{s}_l]$ gives actually rise to an action of the Weyl group $W_n$ on $\mathcal{B}[\bd{s}_l]$, but we do not need this fact in the 
sequel. The next proposition gives a simple expression for $\sigma_i(\bd{\lambda}_l)$ when $\bd{\lambda}_l$ is located at the head of an $i$-chain 
in $\mathcal{B}[\bd{s}_l]$. 

\begin{prop} \label{bij_scopes} Let $\bd{s}_l \in \mZ^l(s)$, $w \in \mathcal{P}(\Fq)$ and $i \in \interv{0}{n-1}$ be such that $w + \alpha_i$ is not 
a weight of $\Fq$. Let $\bd{\lambda}_l \in \Pi^l(\bd{s}_l;w)$ and $\bd{\mu}_l:=\sigma_i(\bd{\lambda}_l)$. Then

\begin{itemize}
\item[\rm (i)] $\bd{\mu}_l$ is the multi-partition obtained by adding to $\bd{\lambda}_l$ all its addable $i$-nodes,
and there are $k_i=(w,\alpha_i)$ of them.
\item[\rm (ii)] We have $|\bd{\mu}_l,\bd{s}_l \rangle = f_i^{(k_i)}. |\bd{\lambda}_l,\bd{s}_l \rangle$ and 
$|\bd{\lambda}_l,\bd{s}_l \rangle = e_i^{(k_i)}. |\bd{\mu}_l,\bd{s}_l \rangle$.
\end{itemize}
\end{prop} 

\begin{proof} See \cite[Prop. 3.5]{Y2}. \end{proof} %réf pour la prop dans Y2!!!!!!!!!!!!!!!!!!!!!!!!!

\subsection{Uglov's canonical bases of the Fock spaces}

Following \cite{U}, the space $\Lambda^s$ can be endowed with an involution $\barre$. Instead of recalling the definition of this 
involution, we give its main properties. They turn out to characterize it completely; see Remark \ref{rq_base_Lambda_s_barre_invariante}.

\begin{prop}[\cite{U}] \label{prop_def_barre} There exists an involution $\barre$ of $\Lambda^s$ such that:
\begin{itemize}
\item[\rm (i)] $\barre$ is a $\mQ$-linear map of $\Lambda^s$ such that for all $u \in \Lambda^s$, $k \in \mZ$, we have 
$\overline{q^k\, u}=q^{-k}\,\overline{u}$.
\item[\rm (ii)] \emph{(Unitriangularity property.)} For all $\lambda \in \Pi$, we have 
$$\overline{| \lambda,s \rangle} \in | \lambda,s \rangle + \ds \bigoplus_{\mu \lhd \lambda}{\mZ[q,q^{-1}] \,|\mu,s \rangle},$$ where $\lhd$ stands 
for the dominance ordering on partitions.
\item[\rm (iii)] For all $\lambda \in \Pi$, we have $\wt(\overline{| \lambda,s \rangle})=\wt(|\lambda,s \rangle)$ and 
$\wtpt(\overline{| \lambda,s \rangle})=\wtpt(|\lambda,s \rangle)$.
\item[\rm (iv)] For all $0 \leq i \leq n-1$,  $0 \leq j \leq l-1$, $m<0$, $v \in \Lambda^s$, we have
$$\overline{f_i.v}=f_i.\overline{v}, \quad \overline{\dot{f}_j.v}=\dot{f}_j.\overline{v} \quad \mbox{and} \quad
\overline{B_m.v}=B_m.\overline{v}.$$ \cqfd
\end{itemize} 
\end{prop}

By \cite{U}, this involution can be computed by straightening non-ordered $q$-wedge products.
However, the number of $q$-wedge products of $2$ factors to be straightened is in general too large for practical computations. We 
shall derive in Section \ref{section_algo_bc} another more efficient algorithm that does not require straightening 
$q$-wedge products. \\

By Propositions \ref{prop_def_barre} (iii) and \ref{fock_spaces_weight_subspaces}, the higher-level Fock spaces
are stable under the involution $\barre$. The involution induced on these spaces will still be 
denoted by $\barre$. Let $\bd{s}_l \in \mZ^l(s)$. For $\bd{\mu}_l \in \Pi^l$, write 
\begin{equation}
\overline{|\bd{\mu}_l,\bd{s}_l \rangle} = \ds \sum_{\bd{\lambda}_l \in \Pi^l}
{a_{\bd{\lambda}_l,\bd{\mu}_l ;\, \bd{s}_l}(q)\, |\bd{\lambda}_l,\bd{s}_l \rangle}
\end{equation}
\noindent with $a_{\bd{\lambda}_l,\bd{\mu}_l ;\,\bd{s}_l}(q) \in \mZ[q,q^{-1}]$, and let 
\begin{equation}
A_{\bd{s}_l}(q):=\bigl(a_{\bd{\lambda}_l,\bd{\mu}_l ;\,\bd{s}_l}(q) \bigr)_{\bd{\lambda}_l,\bd{\mu}_l \in \Pi^l}
\end{equation}
denote the matrix of the involution $\barre$ of $\Fq$ with respect to the standard basis. Since the weight subspaces of $\Fq$ are stable under the 
involution $\barre$, (\ref{poidsl}) implies that $a_{\bd{\lambda}_l,\bd{\mu}_l ;\,\bd{s}_l}(q)$ is zero unless $|\bd{\lambda}_l| = |\bd{\mu}_l|$,
where $|\bd{\lambda}_l|$ (\resp $|\bd{\mu}_l|$) denotes the number of boxes contained in the Young diagram of $\bd{\lambda}_l$ 
(\resp $\bd{\mu}_l$). By Proposition \ref{prop_def_barre} (ii), the matrix $A_{\bd{s}_l}(q)$ is unitriangular. As
a consequence, one can define, by a classical argument, canonical bases as follows.

\begin{thm}[\cite{U}] \label{thm_bases_canoniques_Fqsl}
Let $\bd{s}_l \in \mZ^l(s)$. Then there exists a unique basis
$$\set{G^+(\bd{\lambda}_l, \bd{s}_l) \mid \bd{\lambda}_l \in \Pi^l} \qquad \Bigl( \mbox{\resp}
\set{ G^-(\bd{\lambda}_l, \bd{s}_l) \mid \bd{\lambda}_l \in \Pi^l} \Bigr)$$
of $\Fq$ such that:
$$\begin{array}{lll}
\mbox{\rm (i)} & \overline{G^+(\bd{\lambda}_l, \bd{s}_l)}=G^+(\bd{\lambda}_l, \bd{s}_l) & (\mbox{\resp} 
 \overline{G^-(\bd{\lambda}_l, \bd{s}_l)}=G^-(\bd{\lambda}_l, \bd{s}_l) \; \mbox{\emph{),}}  \\ 
\mbox{\rm (ii)} & G^+(\bd{\lambda}_l, \bd{s}_l) \equiv |\bd{\lambda}_l, \bd{s}_l \rangle \bmod \ q\mathcal{L}^{+}[\bd{s}_l] & (\mbox{\resp}
G^-(\bd{\lambda}_l, \bd{s}_l) \equiv |\bd{\lambda}_l, \bd{s}_l \rangle \bmod \ q^{-1}\mathcal{L}^{-}[\bd{s}_l] \; ), \\ 
\end{array}$$
where $\mathcal{L}^{\epsilon}[\bd{s}_l]:= \ds \bigoplus_{\bd{\lambda}_l \in \Pi^l}{\mZ[q^\epsilon]\,|\bd{\lambda}_l, \bd{s}_l \rangle}$ 
$(\epsilon=\pm 1)$. \cqfd
\end{thm}

For $\epsilon=\pm 1$, define entries $\Delta^{\epsilon}_{\bd{\lambda}_l,\bd{\mu}_l ;\,\bd{s}_l}(q) \in \mZ[q,q^{-1}]$ 
($\bd{\lambda}_l$, $\bd{\mu}_l \in \Pi^l$) by
\begin{equation} \label{def_Delta_q}
G^{\epsilon}(\bd{\mu}_l,\bd{s}_l) = \ds \sum_{\bd{\lambda}_l \in \Pi^l}
{\Delta^{\epsilon}_{\bd{\lambda}_l,\bd{\mu}_l ;\,\bd{s}_l}(q) \,|\bd{\lambda}_l, \bd{s}_l \rangle},
\end{equation}
and denote by 
\begin{equation}
\Delta^{\epsilon}_{\bd{s}_l}(q):=\bigl(\Delta^{\epsilon}_{\bd{\lambda}_l,\bd{\mu}_l ;\,\bd{s}_l}(q) \bigr)_{\bd{\lambda}_l,\bd{\mu}_l \in \Pi^l}
\qquad  (\epsilon = \pm 1)
\end{equation}
the transition matrices between the standard and the canonical bases of $\Fq$. \\

We shall give some (parts of the) canonical bases of different Fock spaces in Section \ref{section_tables}. By \cite{U}, the entries of 
$\Delta^{+}_{\bd{s}_l}(q)$ (\resp $\Delta^{-}_{\bd{s}_l}(q)$) are Kazhdan-Lusztig polynomials of parabolic
submodules of affine Hecke algebras of type $A$, so by \cite{KT}, these polynomials are in $\mN[q]$ (\resp $\mN[p]$). Moreover, both canonical 
bases of $\Fq$ are dual to each other with respect to a certain bilinear form, which gives an inversion formula for Kazhdan-Lusztig polynomials; see 
\cite[Thm. 5.15]{U}. By \cite{U}, the basis 
$\set{G^+(\bd{\lambda}_l, \bd{s}_l) \mid \bd{\lambda}_l \in \Pi^l}$ is a lower global crystal basis (in the sense of \cite{Kas1}) of the integrable 
$\Uq$-module $\Fq$. 

\section{Computation of the canonical bases of $\Fq$} \label{section_algo_bc}
\subsection{A $\barre$-invariant basis $\mathcal{B}$ of $\Lambda^s$}
\begin{notation} \label{notation_barre_invariant_basis}
In this section, we use the following notation.
\begin{itemize} 
\item[*] In this article, we always identify the multi-partition $\bd{\lambda}_l=(\lambda^{(1)},\ldots,\lambda^{(l)}) \in \Pi^l$
with its Young diagram $\set{(i,j,b) \in \mN^* \times \mN^* \times \interv{1}{l} \mid 1 \leq j \leq \lambda^{(b)}_i}\,$.
Let $(i(\gamma),j(\gamma),b(\gamma))$
denote the coordinates of the node $\gamma \in \bd{\lambda}_l$. Let 
\begin{equation}
\partial \bd{\lambda}_l:=\set{(i,j,b) \in \bd{\lambda}_l \mid j=\lambda^{(b)}_i}
\end{equation}
denote the \emph{border} of $\bd{\lambda}_l$, that is the vertical strip made of the rightmost nodes of $\bd{\lambda}_l$. 

\item[*] For $\bd{s}_l \in \mZ^l(s)$, let
\begin{equation}
\Mq:=\Uq.|\bd{\emptyset}_l,\bd{s}_l \rangle.
\end{equation}
Let $\Pi^l(\bd{s}_l)^{°}=\Pi^l(\bd{s}_l,n)^{°}$ denote the set of $l$-multi-partitions indexing the crystal graph of the $\Uq$-module $\Mq$.
For $\bd{s}_l \in \mZ^l(s)$ and $w \in \mathcal{P}(\Fq)$, put
\begin{equation}
\Pi^l(\bd{s}_l;w)^{°}=\Pi^l(\bd{s}_l,n;w)^{°}:= \Pi^l(\bd{s}_l;w) \cap \Pi^l(\bd{s}_l)^{°}.
\end{equation}
\item[*] Let 
\begin{equation} \label{eq_notation_Xln}
\mathcal{X}_{l,n}:=\set{(v_1,\ldots,v_l) \in \mZ^{l} \mid n-1 \geq v_1 \geq \cdots \geq v_l \geq 0}.
\end{equation} 
Note that the $\Uqprime$-module $\Mq$ $\big(\bd{s}_l=(s_1,\ldots,s_l) \in \mZ^l(s)\big)$ is isomorphic to $\Mq[\bd{v}_l]$ for a unique 
multi-charge $\bd{v}_l=(v_1,\ldots,v_l) \in \mathcal{X}_{l,n}$ ; more precisely,
this $\bd{v}_l$ is related to $\bd{s}_l$ by $\Lambda_{v_1}+\cdots+\Lambda_{v_l}=\Lambda_{s_1}+\cdots+\Lambda_{s_l}$.
\fini
\end{itemize}
\end{notation}

\subsubsection{A basis of $\Mq$ $(\bd{s}_l \in \mZ^l(s))$}

In \cite{J1}, Jacon gave a monomial basis of $\Mq[\bd{v}_l]$ with $\bd{v}_l \in \mathcal{X}_{l,n}$.
Using this basis, he derived in \cite{J2} an algorithm for computing the canonical bases of $\Mq[\bd{v}_l]$. We recall here his result and give a basis of the same type for $\Mq$ with $\bd{s}_l \in \mZ^l(s)$.\\

Let $\bd{v}_l \in \mathcal{X}_{l,n}$, $\bd{\lambda}_l \in \Pi^l$ and $0 \leq k \leq n-1$. Let $X_k$ be the set
of all removable $k$-nodes $\gamma$ in $\bd{\lambda}_l$ such that for every $(k-1)$-node $\beta$ in $\partial \bd{\lambda}_l$, we
have $j(\gamma) > j(\beta)$. Here in the definition of $i$-nodes ($i=k$ or $k-1$), the residues are taken with respect to the multi-charge 
$\bd{v}_l$. 
\begin{example}
Take $n=4$, $l=2$, $\bd{v}_l=(3,1)$ and $\bd{\lambda}_l=((4,2),(4,1))$. Then we have
$$X_0=\set{(1,4,2)}, \quad X_1=\emptyset, \quad X_2=\set{(1,4,1)}
\quad \mbox{and} \quad X_3=\emptyset.$$ \fini     
\end{example}

Jacon proved the following result.

\begin{prop}[\cite{J1}, Lemmas 4.2 \& 4.3] \label{prop_def_suites_Jacon} 
Let $\bd{v}_l \! \in \! \mathcal{X}_{l,n}$ and $\bd{\lambda}_l=(\lambda^{(1)},\ldots,\lambda^{(l)}) \! \in \! \Pi^l(\bd{v}_l)^{°}$, $\bd{\lambda}_l \neq \bd{\emptyset}_l$. Let $j_{\max}:=\max \set{\lambda_1^{(b)} \mid 1 \leq b \leq l}$ denote the largest part of $\bd{\lambda}_l$
and for $0 \leq k \leq n-1$, let $X_k$ be the set defined above. Then there exists $k \in \interv{0}{n-1}$ such that
$$X_k \cap \set{ \gamma \in \bd{\lambda}_l \mid j(\gamma)=j_{\max}}$$
is nonempty. Moreover, the multi-partition obtained from $\bd{\lambda}_l$ by removing all the nodes in $X_k$ (for $k$ as above) is in
$\Pi^l(\bd{v}_l)^{°}$. \cqfd
\end{prop}

Thanks to this proposition, one can define recursively a Jacon element $F(\bd{\lambda}_l) \in \Uq$.

\begin{definition} \label{def_suites_Jacon}
Let $\bd{v}_l \in \mathcal{X}_{l,n}$ and $\bd{\lambda}_l \in \Pi^l(\bd{v}_l)^{°}$. If $\bd{\lambda}_l$ is the empty
multi-partition, let $F(\bd{\lambda}_l):=1 \in \Uq$. Otherwise, let
$k \in \interv{0}{n-1}$ be the minimal integer given by Proposition \ref{prop_def_suites_Jacon} and let $X_k$ be corresponding set of removable
$k$-nodes. Then $\bd{\mu}_l:=\bd{\lambda}_l \setminus X_k$ is in $\Pi^l(\bd{v}_l)^{°}$, so by induction, we can define
\begin{equation}
F(\bd{\lambda}_l):=f_{k}^{(\sharp X_k)}F(\bd{\mu}_l).
\end{equation}
Note that $F(\bd{\lambda}_l)$ is a product of divided powers of the $f_k$'s. \fini
\end{definition}
  
\begin{example} Take $n=4$, $l=2$, $\bd{v}_l=(3,1) \in \mathcal{X}_{l,n}$ and $\bd{\lambda}_l=((4,2),(4,1))$. One can check, \eg by computing the 
crystal graph of $\Mq[\bd{v}_l]$, that $\bd{\lambda}_l \in \Pi^l(\bd{v}_l)^{°}$. We have 
$$F(\bd{\lambda}_l)=f_{0}^{(1)}f_{2}^{(1)}f_{1}^{(1)}f_{3}^{(2)}f_{0}^{(2)}f_{2}^{(2)}f_{1}^{(1)}f_{3}^{(1)}.$$ \fini
\end{example}

Jacon proved the following result.

\begin{thm}[\cite{J1}, Prop. 4.6] \label{thm_unitr_Jacon}
Let $\bd{v}_l \in \mathcal{X}_{l,n}$ and $\bd{\lambda}_l \in \Pi^l(\bd{v}_l)^{°}$. Then
$$F(\bd{\lambda}_l).|\bd{\emptyset}_l,\bd{v}_l \rangle \in |\bd{\lambda}_l,\bd{v}_l \rangle +
\ds \bigoplus_{\stackrel{\bd{\mu}_l \in \Pi^l}{a(\bd{\mu}_l)>a(\bd{\lambda}_l)}}{\mZ[q,q^{-1}]\,|\bd{\mu}_l,\bd{v}_l \rangle},$$ 
where $a : \Pi^l \rightarrow \mZ$ denotes Lusztig's $a$-value \cite{BK,J1}. \cqfd
\end{thm} 

We do not use this $a$-value in the sequel, but only the unitriangularity property of this theorem, which implies the following.

\begin{cor} \label{base_barre_invariante_pour_Mq}
Let $\bd{s}_l=(s_1,\ldots,s_l) \in \mZ^l(s)$, and $\bd{v}_l=(v_1,\ldots,v_l)$ be the unique element in $\mathcal{X}_{l,n}$ such that
$\Lambda_{v_1}+\cdots+\Lambda_{v_l}=\Lambda_{s_1}+\cdots+\Lambda_{s_l}$. For $\bd{\lambda}_l \in \Pi^l(\bd{v}_l)^{°}$, let
$F(\bd{\lambda}_l)$ be the Jacon element associated to $\bd{\lambda}_l$ and $\bd{v}_l$. Then 
$$\set{F(\bd{\lambda}_l).|\bd{\emptyset}_l,\bd{s}_l \rangle \mid \bd{\lambda}_l \in \Pi^l(\bd{v}_l)^{°}}$$
is a basis of $\Mq$.
\end{cor}

\begin{proof} Fix a weight $w$ of $\Mq[\bd{v}_l]$. Let 
$B_w:=\set{F(\bd{\lambda}_l).|\bd{\emptyset}_l,\bd{v}_l \rangle \mid \bd{\lambda}_l \in \Pi^l(\bd{v}_l;w)^{°}}$. 
By Theorem \ref{thm_unitr_Jacon}, the transition matrix between the standard basis of $\Mq[\bd{v}_l]\langle w \rangle$ and $B_w$
is unitriangular with respect to a certain ordering given by the $a$-value, which shows that $B_w$ is a basis of $\Mq[\bd{v}_l]\langle w \rangle$.  
This proves the corollary if $\bd{s}_l=\bd{v}_l$. Now, by definition of $\bd{v}_l$, there exists an isomorphism 
$\Mq[\bd{v}_l] \stackrel{\sim}{\longrightarrow} \Mq$ of $\Uqprime$-modules that maps 
$u.|\bd{\emptyset}_l,\bd{v}_l \rangle$ on $u.|\bd{\emptyset}_l,\bd{s}_l \rangle$
($u \in \Uqprime$). The result follows.
\end{proof}

\subsubsection{A basis of $\mathcal{H}.|\bd{\emptyset}_l,\bd{s}_l \rangle$ $(\bd{s}_l \in \mZ^l(s))$}

Recall that the operators $B_m$, $m>0$ (\resp $m<0$) pairwise commute. For any partition $\lambda=(\lambda_1,\ldots,\lambda_r)$, put
\begin{equation} \label{eq_def_B_lambda}
B_{\lambda}:=B_{\lambda_1} \cdots B_{\lambda_r} \in \mathrm{End}(\Lambda^s)
\end{equation}
(by convention, $B_{\emptyset}$ is the identity operator), and define in a similar way $B_{-\lambda}$.

\begin{prop} \label{base_barre_invariante_pour_H_vide}
Let $\bd{s}_l \in \mZ^l(s)$. Then $\set{B_{-\lambda}.|\bd{\emptyset}_l,\bd{s}_l \rangle \mid \lambda \in \Pi}$ is a basis of
$\mathcal{H}.|\bd{\emptyset}_l,\bd{s}_l \rangle$.
\end{prop}

\begin{proof} Let $\mathrm{Sym}$ denote the $\mQ(q)$-vector space of symmetric functions and let $\set{p_{\lambda} \mid \lambda \in \Pi}$ be the 
basis formed by the power sums \cite{Mac}. Let us recall the action of $\mathcal{H}$ on $\mathrm{Sym}$. For $k<0$ let $b_k$ denote the multiplication 
in $\mathrm{Sym}$ by $p_{(-k)}$. For $k>0$ and $f=f(p_{(1)},p_{(2)},\ldots) \in \mathrm{Sym}$, put 
$b_k(f):=\gamma_k \frac{\partial f}{\partial p_{(k)}}$, where $\gamma_k \in \mathcal{Z}(\mathcal{H}) \cong \mQ(q)$ is the scalar such that
$[B_k,B_{-k}]=\gamma_k$. One easily checks that there exists a homomorphism of algebras $\mathcal{H} \rightarrow \mathrm{End}(\mathrm{Sym})$ that
maps $B_k$ on $b_k$ $(k \in \mZ^*)$. This makes $\mathrm{Sym}$ into a simple $\mathcal{H}$-module. Note that the vector 
$v:=|\bd{\emptyset}_l,\bd{s}_l \rangle \in \Fq$ is a singular vector for the action of $\mathcal{H}$, that is, $v \neq 0$ and
$B_k.v=0$ for all $k>0$. By a classical result (see \eg \cite[Lemma 9.13 a)]{Kac}), the linear map:
$\mathrm{Sym} \rightarrow \mathcal{H}.v$, $p_{\lambda} \mapsto B_{-\lambda}.v$ $(\lambda \in \Pi)$, is an isomorphism of $\mathcal{H}$-modules. 
Since $\set{p_{\lambda} \mid \lambda \in \Pi}$ is a basis of $\mathrm{Sym}$, the result follows.
\end{proof}

\subsubsection{A $\barre$-invariant basis $\mathcal{B}$ of $\Lambda^s$}

We need a lemma on bimodules, whose proof is left to the reader.

\begin{lemma} \label{lemma_bimodules}
Let $A$ be $\mQ(q)$-algebra with unity, and $U:=\Uqprime$ or $\Upprime$. Let $V$ be an $(A \otimes U)$-module such that 
the actions of $A \cong (A \otimes 1)$ and $U \cong (1 \otimes U)$ commute. Assume that $V$ is an integrable $U$-module and there exists a
highest weight vector (for the action of $U$), denoted by $v_0$, such that $V = (A \otimes U).v_0 = A.U.v_0$. Let $A^{\bullet} \subset A$ and 
$U^{\bullet} \subset U$ be subsets such that $\set{a.v_0 \mid \ a \in A^{\bullet}}$ is a basis of $A.v_0$ and 
$\set{u.v_0 \mid \ u \in U^{\bullet}}$ is a basis of $U.v_0$. Then $\set{a.u.v_0 \mid a \in A^{\bullet},\ u \in U^{\bullet}}$ is a basis of $V$.\cqfd
\end{lemma}

Recall Notation \ref{notation_barre_invariant_basis} and the definition of the Jacon elements in $\Uqprime$. In a similar way, define 
a Jacon element $\dot{F}(\bd{\lambda}_n) \in \Upprime$ for $\bd{\lambda}_n \in \Pi^l(\bd{v}_n,l)^{°}$, $\bd{v}_n \in \mathcal{X}_{n,l}$. 

\begin{lemma} \label{lemme_base_Lambda_s_barre_invariante}
Keep the notation above. Let $\bd{r}_l \in A_{l,n}(s)$. By \cite{U}, there exists a unique $\bd{r}_n \in A_{n,l}(s)$ such that 
$|\bd{\emptyset}_n,\bd{r}_n \rangle^{\bullet} = |\bd{\emptyset}_l,\bd{r}_l \rangle$. Let $\bd{v}_l$ be the
unique multi-charge in $\mathcal{X}_{l,n}$ such that 
$\Uqprime.|\bd{\emptyset}_l,\bd{v}_l \rangle \cong \Uqprime.|\bd{\emptyset}_l,\bd{r}_l \rangle$ and let $\bd{v}_n$ be the
unique multi-charge in $\mathcal{X}_{n,l}$ such that 
$\Upprime.|\bd{\emptyset}_n,\bd{v}_n \rangle^{\bullet} \cong \Upprime.|\bd{\emptyset}_n,\bd{r}_n \rangle^{\bullet}$. Then 
$$\begin{array}{l}
\set{F(\bd{\lambda}_l) B_{-\mu} \dot{F}(\bd{\lambda}_n).|\bd{\emptyset}_n,\bd{r}_n \rangle^{\bullet} 
\mid \bd{\lambda}_l \in \Pi^l(\bd{v}_l,n)^{°},\ \mu \in \Pi, \ \bd{\lambda}_n \in \Pi^n(\bd{v}_n,l)^{°}} \\[3mm]
= \set{\dot{F}(\bd{\lambda}_n) B_{-\mu} F(\bd{\lambda}_l).|\bd{\emptyset}_l,\bd{r}_l \rangle 
\mid \bd{\lambda}_l \in \Pi^l(\bd{v}_l,n)^{°},\ \mu \in \Pi, \ \bd{\lambda}_n \in \Pi^n(\bd{v}_n,l)^{°}} \\[3mm]
= \set{F(\bd{\lambda}_l) \dot{F}(\bd{\lambda}_n) B_{-\mu}.|\bd{\emptyset}_l,\bd{r}_l \rangle 
\mid \bd{\lambda}_l \in \Pi^l(\bd{v}_l,n)^{°},\ \mu \in \Pi, \ \bd{\lambda}_n \in \Pi^n(\bd{v}_n,l)^{°}}
\end{array}$$
is a basis of the vector space 
$$\Uqprime \otimes \mathcal{H} \otimes \Upprime .|\bd{\emptyset}_n,\bd{r}_n \rangle^{\bullet} 
= \Upprime \otimes \mathcal{H} \otimes \Uqprime .|\bd{\emptyset}_l,\bd{r}_l \rangle
= \Uqprime \otimes \Upprime \otimes \mathcal{H} .|\bd{\emptyset}_l,\bd{r}_l \rangle .$$ 
\end{lemma}

\begin{proof} All the equalities come from the fact that the actions of $\Uqprime$, $\Upprime$ and $\mathcal{H}$ 
pairwise commute. By Corollary \ref{base_barre_invariante_pour_Mq} and Proposition \ref{base_barre_invariante_pour_H_vide},
\begin{itemize}
\item $\set{F(\bd{\lambda}_l).|\bd{\emptyset}_l,\bd{r}_l \rangle \mid \bd{\lambda}_l \in \Pi^l(\bd{v}_l,n)^{°} }$ is a basis of 
$\Uqprime.|\bd{\emptyset}_l,\bd{r}_l \rangle$,
\item $\set{B_{-\mu}.|\bd{\emptyset}_l,\bd{r}_l \rangle \mid \mu \in \Pi}$ is a basis of 
$\mathcal{H}.|\bd{\emptyset}_l,\bd{r}_l \rangle$, and 
\item $\set{\dot{F}(\bd{\lambda}_n).|\bd{\emptyset}_n,\bd{r}_n \rangle^{\bullet} \mid \bd{\lambda}_n \in \Pi^n(\bd{v}_n,l)^{°} }$ is a basis of 
$\Upprime.|\bd{\emptyset}_n,\bd{r}_n \rangle^{\bullet}$.
\end{itemize} 
We now conclude with Lemma \ref{lemma_bimodules}.
\end{proof}

We are now ready to state the following:

\begin{thm} \label{thm_base_Lambda_s_barre_invariante}
With notation of Lemma \ref{lemme_base_Lambda_s_barre_invariante}, the set
$$\begin{array}{cl}
\mathcal{B}:= \set{F(\bd{\lambda}_l) B_{-\mu} \dot{F}(\bd{\lambda}_n).|\bd{\emptyset}_n,\bd{r}_n \rangle^{\bullet}
\mid \bd{r}_n \in A_{n,l}(s),\ \bd{\lambda}_l \in \Pi^l(\bd{v}_l,n)^{°},\ \mu \in \Pi, \ \bd{\lambda}_n \in \Pi^n(\bd{v}_n,l)^{°}} \\[3mm]
= \set{\dot{F}(\bd{\lambda}_n) B_{-\mu} F(\bd{\lambda}_l).|\bd{\emptyset}_l,\bd{r}_l \rangle 
\mid \bd{r}_l \in A_{l,n}(s),\ \bd{\lambda}_l \in \Pi^l(\bd{v}_l,n)^{°},\ \mu \in \Pi, \ \bd{\lambda}_n \in \Pi^n(\bd{v}_n,l)^{°}} \\[3mm]
= \set{F(\bd{\lambda}_l) \dot{F}(\bd{\lambda}_n) B_{-\mu}.|\bd{\emptyset}_l,\bd{r}_l \rangle 
\mid \bd{r}_l \in A_{l,n}(s),\ \bd{\lambda}_l \in \Pi^l(\bd{v}_l,n)^{°},\ \mu \in \Pi, \ \bd{\lambda}_n \in \Pi^n(\bd{v}_n,l)^{°}}
\end{array}$$
is a basis of $\Lambda^s$ that is $\barre$-invariant.
\end{thm}

\begin{proof} The fact that $\mathcal{B}$ is a basis of $\Lambda^s$ comes from Lemma \ref{lemme_base_Lambda_s_barre_invariante} and Theorem 
\ref{thm_dec_Lambda_s}. Let $\bd{r}_l \in \mZ^l(s)$. By (\ref{poidsl}), the subspace of 
$\Fq[\bd{r}_l]$ of weight $\wt(|\bd{\emptyset}_l,\bd{r}_l \rangle)$ is one-dimensional. Therefore, by Proposition 
\ref{prop_def_barre} (ii), the vector $|\bd{\emptyset}_l,\bd{r}_l \rangle$ is $\barre$-invariant. 
Proposition \ref{prop_def_barre} (iv) then implies that every vector in $\mathcal{B}$ is $\barre$-invariant.
\end{proof}

\begin{remark} \label{rq_base_Lambda_s_barre_invariante}
The proof that $\mathcal{B}$ is a basis of $\Lambda^s$ does not use the involution $\barre$. Since every vector in $\mathcal{B}$
is fixed by the involution $\barre$, this determines this involution completely. \fini
\end{remark}     

\subsection{Computation of $\mathcal{B} \cap \Fq \langle v \rangle$ ($\bd{s}_l \in \mZ^l(s)$, $v \in \mathcal{P}(\Fq)$)}
\label{section_computation_LLT_basis}

Let $\bd{s}_l \in \mZ^l(s)$ and $v$ be a weight of $\Fq$. Let $\mathcal{B}$ be the $\barre$-invariant basis of $\Lambda^s$ given by Theorem
\ref{thm_base_Lambda_s_barre_invariante}. The goal of this section is to compute the list of vectors of $\mathcal{B}$ that lie in 
$\Fq \langle v \rangle$. Since the weights of the operators $F(\bd{\lambda}_l)$, $B_{-\mu}$ and $\dot{F}(\bd{\lambda}_n)$ from Theorem
\ref{thm_base_Lambda_s_barre_invariante} are known, it is enough to determine whether the sum of these weights and the weight of 
$|\bd{\emptyset}_l,\bd{r}_l \rangle$ is equal to $v$. This raises no theoretical problem but requires, in view of practical computations, to 
introduce some cumbersome notation. \\

Define a partial ordering on $P:=\bigoplus_{i=0}^{n-1}{\mZ \Lambda_i} \oplus \mZ \delta$ by writing $w' \leq w''$ 
($w'$, $w'' \in P$) if $w''-w' \in \sum_{i=0}^{n-1}{\mN \alpha_i}$.

\begin{definition} \label{def_admissible_weight}
We say that the weight $w \in P$ is \emph{admissible} (with respect to $\Fq \langle v \rangle$) if:
\begin{itemize}
\item[(i)] $w$ is a weight of $\Fq$ such that $v \leq w$.
\item[(ii)] Let $(\bd{r}_n,\dot{w}) \in \mZ^n(s) \times \dot{\mathcal{P}}(\Lambda^s)$ be the unique pair such that 
$\Fq \langle w \rangle = \Fp[\bd{r}_n] \langle \dot{w} \rangle$ (see Proposition \ref{prop_bij_pds_charges}). 
We then require $\bd{r}_n$ to be in $A_{n,l}(s)$. \fini 
\end{itemize}
\end{definition}

\begin{notation} \label{notation_admissible_weight}
Let $w$ be an admissible weight with respect to $\Fq \langle v \rangle$. Then by Definition \ref{def_admissible_weight},
there exist $N_0,\ldots,N_{n-1} \in \mN$ such that $w-v=\sum_{i=0}^{n-1}{N_i\, \alpha_i}$. Let $N(w):=\min_{i}{N_i}$.
The pair $(\bd{r}_n,\dot{w}) \in A_{n,l}(s) \times \dot{\mathcal{P}}(\Lambda^s)$ such that 
$\Fq \langle w \rangle = \Fp[\bd{r}_n] \langle \dot{w} \rangle$ will be denoted by $\big(\bd{r}_n(w),\dot{w}(w)\big)$. Let 
$\bd{r}_l(w) \in A_{l,n}(s)$ be such that $|\bd{\emptyset}_l,\bd{r}_l(w) \rangle = |\bd{\emptyset}_n,\bd{r}_n(w) \rangle^{\bullet}$.
Recall the notation $\mathcal{X}_{l,n}$ and $\mathcal{X}_{n,l}$ from (\ref{eq_notation_Xln}). Then there exists a unique
$\bd{v}_l(w) \in \mathcal{X}_{l,n}$ (\resp $\bd{v}_n(w) \in \mathcal{X}_{n,l}$) such that 
$\Uqprime.|\bd{\emptyset}_l,\bd{v}_l(w) \rangle$ is isomorphic to $\Uqprime.|\bd{\emptyset}_l,\bd{r}_l(w) \rangle$ (\resp
$\Upprime.|\bd{\emptyset}_n,\bd{v}_n(w) \rangle^{\bullet}$ is isomorphic to $\Upprime.|\bd{\emptyset}_n,\bd{r}_n(w) \rangle^{\bullet}$). \fini
\end{notation}

\begin{example} \label{exemple_poids_admissibles}
Take $n=3$, $l=2$, $\bd{s}_l=(3,6)$ and $v = 2 \Lambda_0 - 4 \delta =
\wt(|\bd{\emptyset}_l,\bd{s}_l \rangle) -(\alpha_0+\alpha_1+\alpha_2)$. The following array gives the list of the 
admissible weights $w$ with respect to $\Fq \langle v \rangle$ and the corresponding  
$N(w)$, $\dot{w}(w)$, $\bd{r}_n(w)$, $\bd{v}_n(w)$, $\bd{r}_l(w)$ and $\bd{v}_l(w)$. 
$$\footnotesize \begin{array}{c||c|c|c|c|c|c} 
w  & N(w) & \dot{w}(w) & \bd{r}_n(w) & \bd{v}_n(w) & \bd{r}_l(w) & \bd{v}_l(w) \\
\hline 
2 \Lambda_0 - 3 \delta 		  		  			  & 1 & \stackrel{}{3 \dot{\Lambda}_1 - 3 \dot{\alpha}_1} 
& (3,3,3) & (1,1,1) & (6,3) & (0,0) \\
(2 \Lambda_0 - 3 \delta) -\alpha_0		  			  & 0 & 2 \dot{\Lambda}_0 + \dot{\Lambda}_1 - (\dot{\alpha}_0+3 \dot{\alpha}_1) 
& (4,3,2) & (1,0,0) & (5,4) & (2,1) \\
(2 \Lambda_0 - 3 \delta) -(\alpha_0+\alpha_1+\alpha_2) & 0 & 3 \dot{\Lambda}_1 - (\dot{\alpha}_0+4 \dot{\alpha}_1) 
 & (3,3,3) & (1,1,1) & (6,3) & (0,0) \\
\end{array}$$ \fini
\end{example}

\begin{prop} \label{prop_base_sep_de_Fqsl_barre_invariante}
Let $\bd{s}_l \in \mZ^l(s)$ and $v$ be a weight of $\Fq$. Let $\mathcal{B}$ be the $\barre$-invariant basis of $\Lambda^s$ given by Theorem
\ref{thm_base_Lambda_s_barre_invariante}. Put $\mathbb{B}:=\mathcal{B} \cap \Fq \langle v \rangle$. Then with Notation 
\ref{notation_admissible_weight}, we have
$$\begin{array}{ll}
\mathbb{B} = \ds \bigcup_{\scriptsize w \mbox{ admissible}} \Big\{ \hspace{-3mm} &
F(\bd{\lambda}_l) \dot{F}(\bd{\lambda}_n) B_{-\mu}.|\bd{\emptyset}_l,\bd{r}_l(w) \rangle \ \Big| \ \mu \in \Pi, \, |\mu| \leq N(w), \\
& \bd{\lambda}_n \in \Pi^n(\bd{v}_n(w),l)^{°},\, \wtpt \big(\dot{F}(\bd{\lambda}_n).|\bd{\emptyset}_n,\bd{r}_n(w) \rangle^{\bullet} \big)=\dot{w}(w),
 \\
& \bd{\lambda}_l \in \Pi^l(\bd{v}_l(w),n)^{°},\, \wt \big( F(\bd{\lambda}_l).|\bd{\emptyset}_l,\bd{r}_l(w) \rangle \big)= w'(w,\mu) \Big\},
\end{array}$$
where we put $w'(w,\mu):=v-w+|\mu|\delta+\wt(|\bd{\emptyset}_l,\bd{r}_l(w) \rangle)$.
\end{prop}

\pv \begin{itemize}
\item[*] \underline{Proof of inclusion $\supset$.}
Let $x:=F(\bd{\lambda}_l) \dot{F}(\bd{\lambda}_n) B_{-\mu}.|\bd{\emptyset}_l,\bd{r}_l(w) \rangle$ with 
$w$ admissible and $\bd{\lambda}_l$, $\mu$ and $\bd{\lambda}_n$ as in the right hand-side of the statement of this proposition. We must show that
$x \in \Fq \langle v \rangle$. Put $y:=\dot{F}(\bd{\lambda}_n).|\bd{\emptyset}_n,\bd{r}_n(w) \rangle^{\bullet}$.
Since $\wt \big( F(\bd{\lambda}_l).|\bd{\emptyset}_l,\bd{r}_l(w) \rangle \big)= w'(w,\mu)$
and the operator $B_{-\mu}$ has weight $-|\mu| \delta$, we have
$$\begin{array}{rcl}
\wt(x)
&=& 
\wt \big( F(\bd{\lambda}_l).|\bd{\emptyset}_l,\bd{r}_l(w) \rangle \big) + 
\wt (y) -\wt \big(|\bd{\emptyset}_n,\bd{r}_n(w) \rangle^{\bullet} \big) - |\mu| \delta \\[3mm]
&=&
w'(w,\mu) + \wt(y) -\wt \big(|\bd{\emptyset}_l,\bd{r}_l(w) \rangle \big) - |\mu| \delta \\[3mm]
&=& v - w + \wt(y). 
\end{array}$$
By assumption, we have $\wtpt(y)=\dot{w}(w)$, whence
$y \in \Fp[\bd{r}_n(w)] \langle \dot{w}(w) \rangle =\Fq \langle w \rangle$. This implies 
$\wt(y)=w$ and $\wt(x)=v$. Moreover, since the actions of $\Uqprime$, $\Upprime$ and $\mathcal{H}$ pairwise commute, we have
$$x= F(\bd{\lambda}_l) B_{-\mu}.y \; \in \; \Upprime.\mathcal{H}\,(\Fq \langle w \rangle) \; \subset \; \Fq.$$

\item[*] \underline{Proof of inclusion $\subset$.}
Let $x \in \mathbb{B}$. With notation from Theorem 
\ref{thm_base_Lambda_s_barre_invariante}, let $\bd{r}_n \in A_{n,l}(s)$, $\bd{r}_l \in A_{l,n}(s)$, $\bd{\lambda}_l \in \Pi^l(\bd{v}_l,n)^{°}$, 
$\mu \in \Pi$ and $\bd{\lambda}_n \in \Pi^n(\bd{v}_n,l)^{°}$ be such that 
$$x=F(\bd{\lambda}_l) B_{-\mu} \dot{F}(\bd{\lambda}_n).|\bd{\emptyset}_n,\bd{r}_n \rangle^{\bullet}
   =F(\bd{\lambda}_l) \dot{F}(\bd{\lambda}_n) B_{-\mu} .|\bd{\emptyset}_l,\bd{r}_l \rangle.$$
Since $x \neq 0$, the vector 
$y:=\dot{F}(\bd{\lambda}_n).|\bd{\emptyset}_n,\bd{r}_n \rangle^{\bullet}=\dot{F}(\bd{\lambda}_n).|\bd{\emptyset}_l,\bd{r}_l \rangle$ is nonzero,
so it is a weight vector of the $\Uq$-module $\Lambda^s$. Let us now show that $w:=\wt(y)$ is an admissible weight.
Let $\bd{t}_l \in \mZ^l(s)$ be such that $y \in \Fq[\bd{t}_l]$. Then we have  
$$x=F(\bd{\lambda}_l) B_{-\mu}.y \in \Fq[\bd{t}_l] \big \langle \, w - |\mu| \delta - \textstyle \sum_{i=0}^{n-1}n_i\, \alpha_i \big \rangle,$$
where for $0 \leq i \leq n-1$ we put $n_i:=N_i(\bd{\lambda}_l;\bd{v}_l;n)$. Moreover, since $x \in \Fq \langle v \rangle$ and the spaces
$\Fq[\bd{a}_l]$, $\bd{a}_l \in \mZ^l(s)$ are in direct sum, we must have $\bd{t}_l=\bd{s}_l$ and
$$(*) : \quad v = w - |\mu| \delta - \sum_{i=0}^{n-1}n_i\, \alpha_i  = w - \sum_{i=0}^{n-1}(n_i+|\mu|)\, \alpha_i.$$ 
In particular, we have $y \in \Fq \langle w \rangle$ and $v \leq w$. By Proposition \ref{prop_bij_pds_charges}, there exists a unique pair 
$(\bd{s}_n,\dot{w}) \in \mZ^n(s) \times \dot{\mathcal{P}}(\Lambda^s)$ such that $\Fp \langle \dot{w} \rangle = \Fq \langle w \rangle$. We thus have
$$y = \dot{F}(\bd{\lambda}_n).|\bd{\emptyset}_n,\bd{r}_n \rangle^{\bullet} \in \Fp[\bd{r}_n] \cap \Fq \langle w \rangle
= \Fp[\bd{r}_n] \cap \Fp \langle \dot{w} \rangle \subset \Fp[\bd{r}_n] \cap \Fp.$$ 
Since the $\Fp[\bd{a}_n]$, $\bd{a}_n \in \mZ^n(s)$ are in direct sum, we must have $\bd{s}_n=\bd{r}_n \in A_{n,l}(s)$. As a consequence, $w$ is
an admissible weight. One easily checks that $\dot{w}=\dot{w}(w)$, $\bd{r}_n=\bd{r}_n(w)$, $\bd{v}_n=\bd{v}_n(w)$,
$\bd{r}_l=\bd{r}_l(w)$ and $\bd{v}_l=\bd{v}_l(w)$. By $(*)$ we have $w-v=\sum_{i=0}^{n-1}N_i\,\alpha_i$ with $N_i:=n_i+|\mu|$ $(0 \leq i \leq n-1)$.
Since the $n_i$'s are nonnegative integers, we must have $|\mu|=N_i-n_i \leq N_i$ for all $0 \leq i \leq n-1$, whence 
$|\mu| \leq \min_i(N_i) = N(w)$. Since $y \in \Fq \langle w \rangle = \Fp[\bd{r}_n(w)] \langle \dot{w}(w) \rangle$, we have
$\wtpt \big(\dot{F}(\bd{\lambda}_n).|\bd{\emptyset}_n,\bd{r}_n(w) \rangle^{\bullet} \big)=\dot{w}(w)$. Finally, we have
$$v=\wt(x)= \wt \big( F(\bd{\lambda}_l).|\bd{\emptyset}_l,\bd{r}_l(w) \rangle \big) + 
\wt(y) -\wt \big(|\bd{\emptyset}_n,\bd{r}_n(w) \rangle^{\bullet} \big) - |\mu| \delta$$
and $\wt(y)=w$, which implies $\wt \big( F(\bd{\lambda}_l).|\bd{\emptyset}_l,\bd{r}_l(w) \rangle \big)= w'(w,\mu)$. 
\cqfd
\end{itemize}

Let us now explain how we can apply Proposition \ref{prop_base_sep_de_Fqsl_barre_invariante}. Let $w$ be an admissible weight with 
respect to $\Fq \langle v \rangle$. We compute the list of the multi-partitions 
$\bd{\lambda}_n \in \Pi^n(\bd{v}_n(w),l)^{°}$ such that 
$\wtpt \big(\dot{F}(\bd{\lambda}_n).|\bd{\emptyset}_n,\bd{r}_n(w) \rangle^{\bullet} \big)=\dot{w}(w)$
as follows. The proof of Proposition \ref{prop_base_sep_de_Fqsl_barre_invariante} shows that 
$y:= \dot{F}(\bd{\lambda}_n).|\bd{\emptyset}_n,\bd{r}_n(w) \rangle^{\bullet} \in \Fp[\bd{r}_n(w)] \langle \dot{w}(w) \rangle$ and $y \neq 0$, so 
$\dot{w}(w)$ is a weight of $\Fp[\bd{r}_n(w)]$. By (\ref{poidsptn}), there exist 
$N_0,\ldots,N_{l-1} \in \mN$ such that 
\begin{equation}
\dot{w}(w) = \wtpt(|\bd{\emptyset}_n,\bd{r}_n(w) \rangle^{\bullet}) - \ds \sum_{i=0}^{l-1}{N_i \, \dot{\alpha}_i}.
\end{equation}
Moreover, by definition of the Jacon elements and $(\ref{poidsptn})$, the weight of 
$\dot{F}(\bd{\lambda}_n).|\bd{\emptyset}_n,\bd{r}_n(w) \rangle^{\bullet}$
with $\bd{\lambda}_n \in \Pi^n(\bd{v}_n(w),l)^{°}$ is equal to 
$\wtpt(|\bd{\emptyset}_n,\bd{r}_n(w) \rangle^{\bullet}) - \ds \sum_{i=0}^{l-1}{N_i(\bd{\lambda}_n;\bd{v}_n(w);l)\, \dot{\alpha}_i}$. 
As a consequence, the multi-partitions $\bd{\lambda}_n$ we are looking for are the vertices of the crystal graph of 
$\Up.|\bd{\emptyset}_n,\bd{v}_n(w) \rangle^{\bullet}$ that have exactly $N_i$ $i$-nodes for all $0 \leq i \leq l-1$. In a similar way, we compute the 
list of the multi-partitions $\bd{\lambda}_l \in \Pi^l(\bd{v}_l(w),n)^{°}$ satisfying the condition
$\wt \big( F(\bd{\lambda}_l).|\bd{\emptyset}_l,\bd{r}_l(w) \rangle \big)= w'(w,\mu)$ ($\mu \in \Pi$, $|\mu| \leq N(w)$). \\

We now give an example for Proposition \ref{prop_base_sep_de_Fqsl_barre_invariante}.  

\begin{example}

Take $n=3$, $l=2$, $\bd{s}_l=(3,6)$ and $v = 2 \Lambda_0 - 4 \delta =
\wt(|\bd{\emptyset}_l,\bd{s}_l \rangle) -(\alpha_0+\alpha_1+\alpha_2)$. We compute the basis $\mathbb{B}$ 
of $\Fq \langle v \rangle$ given by Proposition \ref{prop_base_sep_de_Fqsl_barre_invariante}. By Example \ref{exemple_poids_admissibles},
the admissible weights are $2 \Lambda_0-3 \delta$, $(2 \Lambda_0-3 \delta) -\alpha_0$ and 
$(2 \Lambda_0-3 \delta) -(\alpha_0+\alpha_1+\alpha_2)$.

\begin{itemize}
\item[*] \underline{Contribution of $w := 2 \Lambda_0-3 \delta$.} Example \ref{exemple_poids_admissibles} gives 
$\dot{w}(w)=\wtpt(|\bd{\emptyset}_n,\bd{r}_n(w) \rangle^{\bullet}) - 3 \dot{\alpha}_1$, $\bd{v}_n(w)=(1,1,1)$, $\bd{r}_l(w)=(6,3)$ and 
$\bd{v}_l(w)=(0,0)$. By computing the first $4$ layers of the crystal graph of $\Up.|\bd{\emptyset}_n,\bd{v}_n(w) \rangle^{\bullet}$,
we see that the only multi-partition $\bd{\lambda}_n$ in $\Pi^n(\bd{v}_n(w),l)^{°}$ such that 
$\wtpt \big(\dot{F}(\bd{\lambda}_n).|\bd{\emptyset}_n,\bd{r}_n(w) \rangle^{\bullet} \big)=\dot{w}(w)$ is 
$\bd{\lambda}_n:=\bigl((1),(1),(1) \bigr)$. The 
corresponding Jacon element is $$\dot{F}(\bd{\lambda}_n)=\dot{f}_1^{(3)}.$$
Moreover, since $N(w)=1$, the only contributing partitions $\mu$ are $\mu=\emptyset$ and $\mu=(1)$. \\
Take $\mu=\emptyset$. With notation from Proposition \ref{prop_base_sep_de_Fqsl_barre_invariante}, we have
$$w'(w,\mu)=\wt(|\bd{\emptyset}_l,\bd{r}_l(w) \rangle)-(\alpha_0+\alpha_1+\alpha_2),$$ 
so we compute the first $3$ layers of the crystal graph of $\Uq.|\bd{\emptyset}_l,\bd{r}_l(w) \rangle$. 
Doing this, we see that the only multi-partitions $\bd{\lambda}_l \in \Pi^l(\bd{v}_l(w),n)^{°}$ satisfying 
the condition $\wt \big(F(\bd{\lambda}_l).|\bd{\emptyset}_l,\bd{r}_l(w) \rangle \big)=w'(w,\mu)$ are
$\bigl( \emptyset,(3) \bigr)$ and $\bigl( \emptyset,(2,1) \bigr)$. The corresponding
Jacon elements are $f_2f_1f_0$ and $f_1f_2f_0$. We therefore get $2$ vectors
of $\mathbb{B}$, namely 
$$v_1:= \dot{f}_1^{(3)} f_2f_1f_0.|\bd{\emptyset}_l,(6,3) \rangle \qquad \mbox{and} \qquad 
  v_2:= \dot{f}_1^{(3)} f_1f_2f_0.|\bd{\emptyset}_l,(6,3) \rangle.$$
Taking $w=2 \Lambda_0-3 \delta$, $\mu=(1)$ gives another vector of $\mathbb{B}$:
$$v_3:= \dot{f}_1^{(3)} B_{-1}.|\bd{\emptyset}_l,(6,3) \rangle.$$ 
\item[*] \underline{Contribution of the other admissible weights.} In the same way, we get $3$ other vectors of 
$\mathbb{B}$: 
$$v_4:= \dot{f}_1^{(3)} \dot{f}_0 f_2f_1.|\bd{\emptyset}_l,(5,4) \rangle, \quad 
  v_5:= \dot{f}_1^{(3)} \dot{f}_0 f_1f_2.|\bd{\emptyset}_l,(5,4) \rangle \quad \mbox{and} \quad
  v_6:= \dot{f}_1^{(3)} \dot{f}_0\dot{f}_1.|\bd{\emptyset}_l,(6,3) \rangle.$$
\end{itemize}
By Proposition \ref{prop_base_sep_de_Fqsl_barre_invariante}, we therefore have 
$$\begin{array}{rclll}
\mathbb{B} &=& \set{v_1,\ldots,v_6} \\[3mm]
&=& \Big\{  \dot{f}_1^{(3)} f_2f_1f_0.|\bd{\emptyset}_l,(6,3) \rangle,& \dot{f}_1^{(3)} f_1f_2f_0.|\bd{\emptyset}_l,(6,3) \rangle,
	  & \dot{f}_1^{(3)} B_{-1} .|\bd{\emptyset}_l,(6,3) \rangle,\\
&& \quad \dot{f}_1^{(3)} \dot{f}_0 f_2f_1 .|\bd{\emptyset}_l,(5,4) \rangle, &
 \dot{f}_1^{(3)} \dot{f}_0 f_1f_2 .|\bd{\emptyset}_l,(5,4) \rangle, & \dot{f}_1^{(3)} \dot{f}_0\dot{f}_1.|\bd{\emptyset}_l,(6,3) \rangle \Big\}. \\
\end{array}$$ \fini
\end{example}

\subsection{Action of the Heisenberg algebra $\mathcal{H}$ on the vectors $| \bd{\emptyset}_l, \bd{s}_l \rangle$ 
($\bd{s}_l \in \mZ^l(s)$)} \label{section_action_H}

In this section, we explain how to compute the action of $\mathcal{H}$ on the vectors $| \bd{\emptyset}_l, \bd{s}_l \rangle$
($\bd{s}_l \in \mZ^l(s)$) without using the straightening relations. Using the action of $\Upprime$, we can restrict ourselves to
the case where $\bd{s}_l=(s_1,\ldots,s_l)$ is dominant, that is $s_1 \gg \cdots \gg s_l$ 
(see Definition \ref{def_dominant}). By a theorem of Uglov (see Theorem \ref{action_bosons_cas_dominant}),
this case can in turn be reduced to the level-one case. When $l=1$, Leclerc and Thibon \cite{LT2} have given a 
combinatorial formula for the action of $\mathcal{H}$, which avoids the straightening of $q$-wedge products in this case.

\subsubsection{The case $l=1$} 
We now recall the result of Leclerc and Thibon mentioned above. To this aim, we introduce the horizontal ribbon strips. 

\begin{definition} [\cite{LT2}, \S 4.1] Let $\theta$ be a skew diagram. Denote by $\theta \! \! \downarrow$ the horizontal strip made of the bottom 
nodes of the columns of $\theta$. We say that $\theta$ is a \emph{horizontal $N$-ribbon strip of weight $k$} if it can
be tiled by $k$ ribbons of length $N$ whose origins lie in $\theta \! \! \downarrow$. One can check that if such a tiling exists, it is unique.
In this case, the \emph{spin} of $\theta$, denoted by $\mathrm{spin}(\theta)$, is the sum of the heights of the ribbons in the tiling of 
$\theta$. \fini
\end{definition}

Still following \cite{LT2}, we define operators of $\Lambda^s$ as follows. Let $\lambda$, $\mu \in \Pi$ and $N$, $m \in \mN^*$. If 
$\lambda \subset \mu$ and $\mu \setminus \lambda$ is a horizontal $N$-ribbon strip of weight $m$, put 
$L_{\lambda,\mu,m}^{(N)}(q):=(-q)^{-\mathrm{spin}(\theta)}.$ Otherwise, put $L_{\lambda,\mu,m}^{(N)}(q):=0$. Now define $\mathcal{U}_k$,
$\mathcal{V}_k \in \mathrm{End}(\Lambda^s)$ $(k \in \mN^*)$ by
\begin{equation}
\mathcal{U}_k . |\nu,s \rangle = \sum_{\mu \in \Pi}{L_{\mu,\nu,k}^{(n)}(q) \, |\mu,s \rangle} \quad \mbox{and} \quad
\mathcal{V}_k . |\nu,s \rangle = \sum_{\lambda \in \Pi}{L_{\nu,\lambda,k}^{(n)}(q) \, |\lambda,s \rangle} \qquad (\nu \in \Pi).
\end{equation}
Before stating the result of Leclerc and Thibon, let us give some extra notation.
Let $\lambda \in \Pi$. Recall the notation $B_{\lambda}$ and $B_{-\lambda}$ from (\ref{eq_def_B_lambda}), and define in a similar way  
$\mathcal{U}_{\lambda},$  $\mathcal{V}_{\lambda} \in \mathrm{End}(\Lambda^s)$.
Let $\set{h_{\lambda} \mid \lambda \in \Pi}$ denote the set of 
complete symmetric functions and $\set{p_{\lambda} \mid \lambda \in \Pi}$ denote the set of power sums \cite{Mac}. Either set is a basis of the 
space of symmetric functions. Let $(\alpha_{\lambda,\mu})_{\lambda,\mu \in \Pi}$ denote the transition matrix between both bases, namely the matrix 
defined by
\begin{equation} 
p_{\mu} = \sum_{\lambda \in \Pi}{\alpha_{\lambda,\mu}h_{\lambda}} \qquad (\mu \in \Pi).
\end{equation}
The entries $\alpha_{\lambda,\mu}$ can be computed recursively, \eg using Newton formula (see \cite[(2.11)]{Mac}):
\begin{equation} \label{eq_newton_formula}
p_{(k)} = k h_{(k)} - \sum_{i=1}^{k-1}{p_{(i)} h_{(k-i)}} \qquad (k \in \mN^{*}).
\end{equation}
 
\begin{thm}[\cite{LT2}, Thm. 6.4] \label{action_bosons_niveau_1}
 Assume that $l=1$. Then with the notation above, we have
$$B_k = \left\{
\begin{array}{lc}
\ds \sum_{\lambda \in \Pi}{\alpha_{\lambda,(k)} \,\mathcal{U}_{\lambda}} & \mbox{if \quad $k>0$},\\[5mm]
\ds \sum_{\lambda \in \Pi}{\alpha_{\lambda,(|k|)} \,\mathcal{V}_{\lambda}} & \mbox{if \quad $k<0$.}
\end{array}
\right.$$
\cqfd
\end{thm} 

\begin{example} \label{exemple_calcul_bosons_niveau_1}
Take $l=1$, $n=2$, $k=-2$ and $s \in \mZ$. By Theorem \ref{action_bosons_niveau_1} and (\ref{eq_newton_formula}), we have
$B_{-2}=2 \mathcal{V}_2 - \mathcal{V}_1^2$. One computes

$$\footnotesize \begin{array}{rcl}
\mathcal{V}_1.| \emptyset,s \rangle &=& | (2),s \rangle - q^{-1} | (1,1),s \rangle, \\[3mm]
\mathcal{V}_1.| (2),s \rangle &=& | (4),s \rangle + | (2,2),s \rangle -q^{-1} | (2,1,1),s \rangle, \\[3mm]
\mathcal{V}_1.| (1,1),s \rangle &=&  | (3,1),s \rangle -q^{-1} | (2,2),s \rangle - q^{-1}| (1,1,1,1),s \rangle, \\[3mm]
\mathcal{V}_1^2.| \emptyset,s \rangle &=& | (4),s \rangle -q^{-1} | (3,1),s \rangle + (1+q^{-2}) | (2,2),s \rangle 
-q^{-1} | (2,1,1),s \rangle + q^{-2} | (1,1,1,1),s \rangle, \\[3mm]
\mathcal{V}_2.| \emptyset,s \rangle &=& | (4),s \rangle -q^{-1} | (3,1),s \rangle + q^{-2} | (2,2),s \rangle, \\
\end{array}$$
whence $$B_{-2}.| \emptyset,s \rangle = | (4),s \rangle -q^{-1} | (3,1),s \rangle + (-1+q^{-2}) | (2,2),s \rangle + q^{-1} | (2,1,1),s \rangle 
- q^{-2}| (1,1,1,1),s \rangle,$$
as one can check using the straightening relations. \fini
\end{example}

\subsubsection{The dominant case}

We shall now see that by a result of \cite{U} (see Theorem \ref{action_bosons_cas_dominant}), it is possible to compute the action of $B_m$ 
($m \in \mZ^{*}$) on certain vectors of the standard basis of $\Lambda^s$ by reducing to the case $l=1$. In this case, we can apply subsequently 
Theorem \ref{action_bosons_niveau_1} and therefore avoid using the straightening relations. Before stating Theorem \ref{action_bosons_cas_dominant},
let us introduce the following definition and notation.

\begin{definition}[\cite{U}] \label{def_dominant}
 Let $M \in \mN^{*}$, $\bd{\lambda}_l \in \Pi^l$ and $\bd{s}_l=(s_1,\ldots,s_l) \in \mZ^l$. We say that the pair
$(\bd{\lambda}_l,\bd{s}_l)$ is \emph{$M$-dominant} if for all $1 \leq i \leq l-1$, we have
$$s_i-s_{i+1} \geq M + |\bd{\lambda}_l|.$$
\fini
\end{definition}

Recall that $\Lambda^s[n,1]$ denotes the space of $q$-wedge products of charge $s$ in which we take $l=1$. If $x$ is an operator of
$\Lambda^s[n,1]$ that acts on the standard basis by 
\begin{equation}
x.|\lambda,s \rangle = \ds \sum_{\mu \in \Pi}{x_{\lambda,\mu}^{(s)} \, |\mu,s \rangle} \qquad 
\big(\lambda \in \Pi, \, x_{\lambda,\mu}^{(s)} \in \mQ(q)\big),
\end{equation}
define for $b \in \interv{1}{l}$ an operator $x[b]$ of $\Lambda^s[n,l]$ by 
\begin{equation}
\begin{array}{l}
x[b]. | \bd{\lambda}_l, \bd{s}_l \rangle = \ds \sum_{\mu \in \Pi}{x_{\lambda^{(b)},\mu}^{(s_b)} \,
\left | (\lambda^{(1)},\ldots,\lambda^{(b-1)},\mu,\lambda^{(b+1)},\ldots,\lambda^{(l)}), \bd{s}_l \right \rangle} \\[5mm]
\qquad (\bd{\lambda}_l=(\lambda^{(1)},\ldots,\lambda^{(l)}) \in \Pi^l, \quad \bd{s}_l=(s_1,\ldots,s_l) \in \mZ^l(s)) \,;
\end{array}
\end{equation}
in other words, $x[b]$ acts as $x$ on the $b$-th component and trivially on the other components.
  
\begin{thm}[\cite{U}, Prop 5.3] \label{action_bosons_cas_dominant}
Let $m \in \mZ^{*}$. Assume that $(\bd{\lambda}_l, \bd{s}_l)$ is $n\widetilde{m}$-dominant, where $\widetilde{m}:=\max(0,-m)$. Then we have 
$$B_m.|\bd{\lambda}_l, \bd{s}_l \rangle = \sum_{b=1}^l {q^{(b-1)|m|} \,B_m[b].|\bd{\lambda}_l, \bd{s}_l \rangle}.$$ \cqfd
\end{thm} 

\begin{example} 

Take $n=2$, $l=2$, $m=-2$, $\bd{\lambda}_l=\bd{\emptyset}_l$ and $\bd{s}_l=(2,-2)$. By Theorem 
\ref{action_bosons_cas_dominant} and Example \ref{exemple_calcul_bosons_niveau_1}, we have

$$\begin{array}{rcl}
B_{-2}.|\bd{\emptyset}_l,\bd{s}_l \rangle 
&=& B_{-2}[1].|\bd{\emptyset}_l,\bd{s}_l \rangle + q^2 B_{-2}[2].|\bd{\emptyset}_l,\bd{s}_l \rangle \\[3mm]
 &=& |((4), \emptyset),\bd{s}_l \rangle -q^{-1} |((3,1), \emptyset),\bd{s}_l \rangle + 
 (-1+q^{-2}) |((2,2), \emptyset),\bd{s}_l \rangle \\[2mm]
 && + q^{-1} |((2,1,1), \emptyset),\bd{s}_l \rangle - q^{-2} |((1,1,1,1), \emptyset),\bd{s}_l \rangle \\[2mm]
 &&+ q^2|(\emptyset,(4)),\bd{s}_l \rangle -q |(\emptyset,(3,1)),\bd{s}_l \rangle + 
(-q^2+1) |(\emptyset,(2,2)),\bd{s}_l \rangle \\[2mm]
&&+ q |(\emptyset,(2,1,1)),\bd{s}_l \rangle - |(\emptyset,(1,1,1,1)),\bd{s}_l \rangle,
\end{array}$$
as one can check using the straightening relations. \fini
\end{example}

\subsubsection{Action of $\mathcal{H}$ on $| \bd{\emptyset}_l, \bd{s}_l \rangle$, $\bd{s}_l \in \mZ^l(s)$}

Let $\lambda=(\lambda_1,\ldots,\lambda_r) \in \Pi$ and $\bd{s}_l \in \mZ^l(s)$. Recall that we put  
$B_{-\lambda}:=B_{-\lambda_1} \cdots B_{-\lambda_r} \in \mathcal{H}$. We derive here a method for computing 
$B_{-\lambda}.|\bd{\emptyset}_l, \bd{s}_l \rangle$ without using the straightening relations. We proceed in four steps. 
The first two steps use some results that we prove in this section.

\begin{enumerate} 
\item We find $\bd{t}_l \in \dot{W}_l.\bd{s}_l$ such that $\bd{t}_l$ is $n|\lambda|$-dominant; see Lemma 
\ref{lemme_existence_multicharge_dominante}.
\item We find $\dot{u} \in \Upprime$ (more precisely, $\dot{u}$ is a monomial in the Chevalley generators of $\Upprime$) such that
$|\bd{\emptyset}_l, \bd{s}_l \rangle = \dot{u}.| \bd{\emptyset}_l, \bd{t}_l \rangle$\,; see Proposition \ref{prop_dec_exterieure_multipartitions}.
\item We claim that we can compute
$B_{-\lambda}.| \bd{\emptyset}_l, \bd{t}_l \rangle$ by repeated applications of Theorem \ref{action_bosons_cas_dominant}. Assume
that $r=2$, \ie $\lambda=(\lambda_1,\lambda_2)$ has only $2$ parts (the general case follows by induction on $r$). By Theorem 
\ref{action_bosons_cas_dominant}, $B_{-\lambda_1}.| \bd{\emptyset}_l, \bd{t}_l \rangle$ is a linear combination of vectors of the form
$| \bd{\mu}_l, \bd{t}_l \rangle$ with $\bd{\mu}_l \in \Pi^l$. Let $\bd{\mu}_l \in \Pi^l$ be such that $| \bd{\mu}_l, \bd{t}_l \rangle$ appears
in this linear combination. Put $w:=\wt(|\bd{\emptyset}_l, \bd{t}_l \rangle)$. Since $B_{-\lambda_1}.\Lambda^s \langle w \rangle \subset
\Lambda^s \langle w-\lambda_1 \delta \rangle$, Equation (\ref{poidsl}) implies that 
$|\bd{\mu}_l|=n\lambda_1$, so $(\bd{\mu}_l, \bd{t}_l)$ is $n\lambda_2$-dominant. We can therefore apply 
Theorem \ref{action_bosons_cas_dominant} to compute $B_{-\lambda_2}.|\bd{\mu}_l, \bd{t}_l \rangle$, which proves the claim.
\item Since the actions of $\Upprime$ and $\mathcal{H}$ commute, we have 
$B_{-\lambda}.| \bd{\emptyset}_l, \bd{s}_l \rangle = \dot{u}.(B_{-\lambda}.| \bd{\emptyset}_l, \bd{t}_l \rangle)$.  
Finally, we compute the action of $\dot{u}$ on $B_{-\lambda}.| \bd{\emptyset}_l, \bd{t}_l \rangle$ by using the indexation of
the standard basis of $\Lambda^s$ by $n$-multi-partitions.  
\end{enumerate}

We now state and prove Lemma \ref{lemme_existence_multicharge_dominante} and Proposition \ref{prop_dec_exterieure_multipartitions}. We give
an example of application at the end of this section.

\begin{lemma} \label{lemme_existence_multicharge_dominante}
Let $\bd{s}_l \in \mZ^l(s)$ and $M \in \mN^{*}$. Then there exists an $nM$-dominant multi-charge $\bd{t}_l$ that is $\dot{W}_l$-conjugated to 
$\bd{s}_l$. 
\end{lemma}

\begin{proof} Let $\lambda:=((l-1)M,(l-2)M,\ldots,M,-l(l-1)M/2) \in \mZ^l$ and $\dot{\tau} \in \mathrm{End}(\mQ^l)$ be the translation by 
$n \lambda$. Since $\lambda \in \mZ^l(0)$, we have $\dot{\tau} \in \dot{W}_l$. We can choose $\dot{\sigma} \in \mathfrak{S}_l \subset \dot{W}_l$ such 
that $\bd{a}_l=(a_1,\ldots,a_l):=\dot{\sigma}.\bd{s}_l$ satisfies $a_1 \geq a_2 \geq \ldots \geq a_l$. Then 
$\bd{t}_l:=\dot{\tau} \dot{\sigma}.\bd{s}_l$ is $nM$-dominant and $\dot{W}_l$-conjugated to $\bd{s}_l$.
\end{proof}

Recall the definition of $\sigma_i(\bd{\lambda}_l) \in \Pi^l$ ($\bd{\lambda}_l \in \Pi^l$, $0 \leq i \leq n-1$) from Definition
\ref{def_sigma_i}. In a similar way, one defines for $0 \leq i \leq l-1$ and $\bd{\lambda}_n \in \Pi^n$ an $n$-multi-partition
$\dot{\sigma}_i(\bd{\lambda}_n)$.

\begin{lemma} \label{lemme_dec_exterieure_multipartitions}
Let $\bd{r}_l \in A_{l,n}(s)$. Let $\bd{r}_n$ be the unique multi-charge in $A_{n,l}(s)$ such that 
$|\bd{\emptyset}_n,\bd{r}_n \rangle^{\bullet} = |\bd{\emptyset}_l,\bd{r}_l \rangle$.
Let $\dot{\sigma}=\dot{\sigma}_{i_r} \ldots \dot{\sigma}_{i_1}$ be a reduced expression of $\dot{\sigma} \in \dot{W}_l$. Define a
sequence $(\bd{\lambda}_n^{(0)},\ldots,\bd{\lambda}_n^{(r)})$ of $n$-multi-partitions by
$$\bd{\lambda}_n^{(0)}:=\bd{\emptyset}_n \qquad \mbox{and} \qquad \bd{\lambda}_n^{(j)}:=\dot{\sigma}_{i_j}(\bd{\lambda}_n^{(j-1)}) 
\quad (1 \leq j \leq r).$$ 
Then for all $1 \leq j \leq r$, $\wtpt \big(| \bd{\lambda}_n^{(j-1)},\bd{r}_n \rangle^{\bullet} \big)+\dot{\alpha}_{i_j}$ is not a weight of 
$\Fp[\bd{r}_n]$, and we have
$$|\bd{\lambda}_n^{(r)},\bd{r}_n \rangle^{\bullet} = | \bd{\emptyset}_l, \dot{\sigma}.\bd{r}_l \rangle.$$ 
\end{lemma}

\begin{proof} We proceed by induction on $r$. If $r=0$, the equality 
$|\bd{\lambda}_n^{(r)},\bd{r}_n \rangle^{\bullet} = | \bd{\emptyset}_l, \dot{\sigma}.\bd{r}_l \rangle$ is obvious and (\ref{poidsptn}) in
Prop. \ref{prop_formulas_weights} shows that for all $i_1 \in \interv{0}{l-1}$,
$\wtpt \big(| \bd{\lambda}_n^{(0)},\bd{r}_n \rangle^{\bullet} \big)+\dot{\alpha}_{i_1}$ is not a weight of 
$\Fp[\bd{r}_n]$. Assume now that $r>0$. Put $\dot{\tau}:=\dot{\sigma}_{i_{r-1}} \cdots \dot{\sigma}_{i_1}$, 
$\bd{s}_l=(s_1,\ldots,s_l):=\dot{\sigma}.\bd{r}_l$, $\bd{t}_l=(t_1,\ldots,t_l):=\dot{\tau}.\bd{r}_l$, 
$u:=|\bd{\lambda}_n^{(r-1)},\bd{r}_n \rangle^{\bullet}$, $\dot{w}:=\wtpt(u)$,
$x:=| \bd{\emptyset}_l, \bd{s}_l \rangle$, $y:=|\bd{\lambda}_n^{(r)},\bd{r}_n \rangle^{\bullet}$ and $i:=i_r$. By induction, 
we have $u=| \bd{\emptyset}_l, \bd{t}_l \rangle$. We must show that $x=y$ and $\dot{w}+\dot{\alpha}_i$ is not a weight of $\Fp[\bd{r}_n]$. \\
We claim that $x$ and $y$ lie in $\Fp[\bd{r}_n]\langle \dot{\sigma}_i.\dot{w} \rangle$. The claim for $y$ follows by the definition of 
$\bd{\lambda}_n^{(r)} = \dot{\sigma}_i(\bd{\lambda}_n^{(r-1)}) \in \Pi^n$. 
Since $\bd{s}_l \in \dot{W}_l.\bd{t}_l$, we have $\Lambda_{s_1}+\cdots+\Lambda_{s_l}=\Lambda_{t_1}+\cdots+\Lambda_{t_l}$. Therefore (\ref{poidsl}) 
implies $\wt(x) \equiv \wt(u) \bmod \mZ \delta$. Moreover, Proposition \ref{fock_spaces_weight_subspaces} implies that
$u = |\bd{\lambda}_n^{(r-1)},\bd{r}_n \rangle^{\bullet}$ is in $\Fp[\bd{r}_n] = \ds 
\bigoplus_{d \in \mZ}{\Lambda^s \langle a_0 \Lambda_0+\cdots+a_{n-1}\Lambda_{n-1}+d \delta \rangle}$ with 
$(a_0,\ldots,a_{n-1}):=\theta_{n,l}(\bd{r_n})$. This implies that $x \in \Fp[\bd{r}_n]$. Now we show easily that
$\wtpt(x)=\dot{\sigma}_i.\wtpt(u)=\wtpt(y)$, whence the claim for $x$. Put 
$\dot{w}_{\emptyset}:=\wtpt\big( |\bd{\emptyset}_n,\bd{r}_n \rangle^{\bullet} \big)$. By (\ref{poidsptn}), the weight subspace 
$\Fp[\bd{r}_n]\langle \dot{w}_{\emptyset} \rangle$ is one-dimensional. Since the formal character of the integrable $\Up$-module $\Fp[\bd{r}_n]$ is
$\dot{W}_l$-invariant, we have 
$$\mathrm{dim}(\Fp[\bd{r}_n] \langle \dot{\sigma}_{i}.\dot{w} \rangle) = 
\mathrm{dim}(\Fp[\bd{r}_n] \langle \dot{\sigma}.\dot{w}_{\emptyset} \rangle) =
\mathrm{dim}(\Fp[\bd{r}_n] \langle \dot{w}_{\emptyset} \rangle)=1.$$
As a consequence, $x$ and $y$ are vectors of the standard basis of a one-dimensional subspace of $\Lambda^s$, whence $x=y$. 

Let us now show that $\dot{w}+\dot{\alpha}_i$ is not a weight of $\Fp[\bd{r}_n]$. Put
$$\dot{\varepsilon}_i:=\max \set{m \in \mN \mid \dot{e}_i^m.x \neq 0} \qquad \mbox{and} \qquad 
\dot{\varphi}_i:=\max \{m \in \mN \mid \dot{f}_i^m.x \neq 0 \}.$$
The theory of $U_p(\gsl{2})$-modules yields 
$$( \dot{\wt}(x),\, \dot{\alpha}_i ) = \dot{\varphi}_i - \dot{\varepsilon}_i.$$
By (\ref{poidsptl}), we have $( \dot{\wt}(x),\, \dot{\alpha}_i ) = s_i-s_{i+1}$, where we put $s_0:=s_l+n$ if $i=0$. 
We thus have $\dot{\varphi}_i - \dot{\varepsilon}_i = s_i-s_{i+1}$. Since $\ell(\dot{\sigma}) = \ell(\dot{\tau})+1$, a classical result 
about Coxeter groups (see \eg \cite[Lemma 2.1 (iii)]{D}) shows that $s_{i+1}>s_i$ (even if $i=0$). As a consequence, 
$\dot{\varepsilon}_i > \dot{\varphi}_i \geq 0$, therefore $\wtpt(x) + \dot{\alpha}_i$ and 
$\dot{\sigma}^{-1}.(\wtpt(x) + \dot{\alpha}_i)= \dot{w}_{\emptyset} + \dot{\sigma}^{-1}.\dot{\alpha}_i$ are weights of $\Fp[\bd{r}_n]$.
By (\ref{poidsptn}), we have $\dot{\sigma}^{-1}.\dot{\alpha}_i \in \sum_{j=0}^{l-1}{\mN \dot{\alpha}_j}$ and 
$\dot{w}_{\emptyset} - \dot{\sigma}^{-1}.\dot{\alpha}_i$ is not a weight of $\Fp[\bd{r}_n]$. Therefore 
$\dot{w}+\dot{\alpha}_i=\dot{\tau}.(\dot{w}_{\emptyset} - \dot{\sigma}^{-1}.\dot{\alpha}_i)$ is not a weight of $\Fp[\bd{r}_n]$ either.
\end{proof}

The following result is a refinement of \cite[Cor. 4.9]{U}.
 
\begin{prop} \label{prop_dec_exterieure_multipartitions} 
Let $\bd{r}_l \in A_{l,n}(s)$ and $\bd{s}_l=(s_1,\ldots,s_l) \in \dot{W}_l. \bd{r}_l$. Then there exist integers $k_1,\ldots,k_r \in \mN$,
$i_1,\ldots,i_r \in \interv{0}{l-1}$ (which can be explicitly calculated) such that 
$$|\bd{\emptyset}_l,\bd{r}_l \rangle = \dot{e}_{i_1}^{(k_1)} \cdots \dot{e}_{i_r}^{(k_r)}.|\bd{\emptyset}_l,\bd{s}_l \rangle \qquad \mbox{and} \qquad
  |\bd{\emptyset}_l,\bd{s}_l \rangle = \dot{f}_{i_r}^{(k_r)} \cdots \dot{f}_{i_1}^{(k_1)}.|\bd{\emptyset}_l,\bd{r}_l \rangle.$$
\end{prop}

\begin{proof}  Let $\dot{\sigma} \in \dot{W}_l$ be the element of minimal length such that $\bd{s}_l =(s_1,\ldots,s_l)= \dot{\sigma}.\bd{r}_l$. 
We argue by induction on the length $r$ of $\dot{\sigma}$. If $r=0$, we have $\bd{s}_l=\bd{r}_l$ and there is nothing to prove. Assume now that $r>0$. We 
compute a reduced expression $$\dot{\sigma}:=\dot{\sigma}_{i_r} \cdots \dot{\sigma}_{i_1}$$ of $\dot{\sigma}$ as follows. Since $r>0$, 
we have $\bd{s}_l \notin A_{l,n}(s)$, so by (\ref{definition_A_l_n_s}) there exists $i_r \in \interv{0}{l-1}$ such that $s_{i_r}<s_{i_r+1}$ 
(here we put, as usual, $s_0:=n+s_l$). Let $\dot{\tau}:=\dot{\sigma}_{i_r} \dot{\sigma}$. By \cite[Lemma 2.1 (iii)]{D}, we have 
$\ell(\dot{\tau})<\ell(\dot{\sigma})$. By induction we get a reduced expression $\dot{\tau}=\dot{\sigma}_{i_{r-1}} \cdots \dot{\sigma}_{i_1}$
of $\dot{\tau}$, so $\dot{\sigma}=\dot{\sigma}_{i_r} \cdots \dot{\sigma}_{i_1}$ is a reduced expression of $\dot{\sigma}$. Now let    
$\bd{r}_n \in A_{n,l}(s)$ be such that $|\bd{\emptyset}_n,\bd{r}_n \rangle^{\bullet} = |\bd{\emptyset}_l,\bd{r}_l \rangle$. 
By Lemma \ref{lemme_dec_exterieure_multipartitions}, we have 
$| \bd{\emptyset}_l,\, \bd{s}_l \rangle = |\dot{\sigma}_{i_r} \ldots \dot{\sigma}_{i_1}(\bd{\emptyset}_n),\bd{r}_n \rangle^{\bullet}$. Let
$$k_j:=|\dot{\sigma}_{i_{j}} \ldots \dot{\sigma}_{i_1}(\bd{\emptyset}_n)|-|\dot{\sigma}_{i_{j-1}} \ldots \dot{\sigma}_{i_1}(\bd{\emptyset}_n)|
\qquad (1 \leq j \leq r).$$
By Proposition \ref{bij_scopes}, we have
$$|\bd{\emptyset}_l,\bd{r}_l \rangle = |\bd{\emptyset}_n,\bd{r}_n \rangle^{\bullet} = \dot{e}_{i_1}^{(k_1)} \cdots \dot{e}_{i_r}^{(k_r)}.
|\dot{\sigma}_{i_r} \ldots \dot{\sigma}_{i_1}(\bd{\emptyset}_n),\bd{r}_n \rangle^{\bullet} 
= \dot{e}_{i_1}^{(k_1)} \cdots \dot{e}_{i_r}^{(k_r)}.|\bd{\emptyset}_l,\bd{s}_l \rangle,$$
and the proof of the other equality is similar.
\end{proof}

\begin{example} Take $s=0$, $n=2$, $l=3$, $\lambda=(1)$ and $\bd{s}_l=(1,-1,0)$. We apply the method derived in this section for computing 
$B_{-\lambda}.| \bd{\emptyset}_l, \bd{s}_l \rangle$. Note that $\bd{s}_l$ is not $n|\lambda|$-dominant, so we cannot directly apply  
Theorem \ref{action_bosons_cas_dominant}. The computation is done in four steps.

\begin{enumerate} 
\item Following the proof of Lemma \ref{lemme_existence_multicharge_dominante}, put $\bd{t}_l:=(3,1,-4)$. Then $\bd{t}_l$ 
is $n|\lambda|$-dominant and $\dot{W}_l$-conjugated to $\bd{s}_l$.
\item As in the proof of Proposition \ref{prop_dec_exterieure_multipartitions}, we compute $\dot{u} \in \Upprime$ such that 
$|\bd{\emptyset}_l, \bd{s}_l \rangle = \dot{u}.| \bd{\emptyset}_l, \bd{t}_l \rangle$. The element in $A_{l,n}(s)$ that is 
$\dot{W}_l$-conjugated to  $\bd{s}_l$ is $\bd{r}_l:=(1,0,-1)$. The following array shows the computation of a reduced expression of $\dot{\sigma}$,
where $\dot{\sigma} \in \dot{W}_l$ is the element of minimal length such that $\bd{s}_l=\dot{\sigma}.\bd{r}_l$. Each line contains a multi-charge
$\bd{a}_l=(a_1,\ldots,a_l)$ and an integer $i \in \interv{0}{l-1}$ such that $a_i<a_{i+1}$ (if it exists). The multi-charge on the next line is 
$\dot{\sigma}_i.\bd{a}_l$.
$$\begin{array}{ll}
(3,1,-4)  & i=0 \\
(-2,1,1)  & i=1 \\
(1,-2,1) & i=2 \\
(1,1,-2) & i=0 \\
(0,1,-1)  & i=1 \\
(1,0,-1)
\end{array}$$
We get $\bd{s}_l = \dot{\sigma}_{i_r} \cdots \dot{\sigma}_{i_1}.\bd{r}_l$ with $r=5$ and $(i_r,\ldots,i_1)=(0,1,2,0,1)$. The element 
$\bd{r}_n \in A_{n,l}(s)$ such that $|\bd{\emptyset}_l,\bd{r}_l \rangle = |\bd{\emptyset}_n,\bd{r}_n \rangle^{\bullet}$ is $\bd{r}_n:=(1,-1)$.
We compute recursively the $n$-multi-partitions $\bd{\lambda}_n^{(j)}:= \dot{\sigma}_{i_j} \cdots \dot{\sigma}_{i_1} (\bd{\emptyset}_n)$
($0 \leq j \leq r$) from $\bd{\lambda}_n^{(0)}=\bd{\emptyset}_n$. The proof of Proposition \ref{prop_dec_exterieure_multipartitions} shows that
for $1 \leq j \leq r$, $k_j$ is the number of addable $i_j$-nodes of $\bd{\lambda}_n^{(j-1)}$ and $\bd{\lambda}_n^{(j)}$ is obtained from 
$\bd{\lambda}_n^{(j-1)}$ by adding these nodes. The computation of the integers $k_j$ $(1 \leq j \leq r)$ is displayed in the following array.

$$\begin{array}{cccc}
j & \bd{\lambda}_n^{(j)} & i_{j+1} & k_{j+1} \\ \hline
0 & ((),()) & 1 & 1 \\
1 & ((1),()) & 0 & 1 \\
2 & ((1,1),()) & 2 & 3 \\
3 & ((2,1,1),(1)) & 1 & 3 \\
4 & ((2,2,1,1),(1,1)) & 0 & 5 \\
5 & ((3,2,2,1,1),(2,1,1)) \\
\end{array}$$ 
By Proposition \ref{prop_dec_exterieure_multipartitions}, we therefore have 
$$|\bd{\emptyset}_l,\bd{r}_l \rangle = 
\dot{e}_{1}^{(1)} \dot{e}_{0}^{(1)} \dot{e}_{2}^{(3)} \dot{e}_{1}^{(3)} \dot{e}_{0}^{(5)}.
|\bd{\emptyset}_l,\bd{t}_l \rangle.$$ 
In the same way we get
$$|\bd{\emptyset}_l,\bd{s}_l \rangle = \dot{f}_{2}^{(1)}.|\bd{\emptyset}_l,\bd{r}_l \rangle.$$ 
As a consequence, we have $|\bd{\emptyset}_l, \bd{s}_l \rangle = \dot{u}.| \bd{\emptyset}_l, \bd{t}_l \rangle$ with
$$\dot{u}=\dot{f}_{2}^{(1)}\dot{e}_{1}^{(1)} \dot{e}_{0}^{(1)} \dot{e}_{2}^{(3)} \dot{e}_{1}^{(3)} \dot{e}_{0}^{(5)}.$$
\item By Theorem \ref{action_bosons_cas_dominant}, we have
$$\begin{array}{rcl}
B_{-\lambda}.|\bd{\emptyset}_l,\bd{t}_l\rangle 
&=& |((2),\emptyset,\emptyset),\bd{t}_l \rangle -q^{-1} |((1,1),\emptyset,\emptyset),\bd{t}_l \rangle \\
&&+\ q|(\emptyset,(2),\emptyset),\bd{t}_l \rangle - |(\emptyset,(1,1),\emptyset),\bd{t}_l \rangle \\
&&+\ q^2|(\emptyset,\emptyset,(2)),\bd{t}_l \rangle -q |(\emptyset,\emptyset,(1,1)),\bd{t}_l \rangle. 
\end{array}$$
\item We only have to apply $\dot{u}$ to get
$$\begin{array}{rcl}
B_{-\lambda}.| \bd{\emptyset}_l, \bd{s}_l \rangle &=& \dot{u}.(B_{-\lambda}.| \bd{\emptyset}_l, \bd{t}_l \rangle) \\[2mm]
&=& |((2),\emptyset,\emptyset),\bd{s}_l \rangle + |(\emptyset,\emptyset,(2)),\bd{s}_l \rangle + (q-q^{-1}) |((1),\emptyset,(1)),\bd{s}_l \rangle \\
&&-\ |(\emptyset,\emptyset,(1,1)),\bd{s}_l \rangle + q^2 |(\emptyset,(2),\emptyset),\bd{s}_l \rangle-|((1,1),\emptyset,\emptyset),\bd{s}_l \rangle\\
&& -\ q |(\emptyset,(1,1),\emptyset),\bd{s}_l \rangle,   
\end{array}$$
as one can check using the straightening relations. \fini 
\end{enumerate}
\end{example}

\subsection{End of the algorithm}

Let $\bd{s}_l \in \mZ^l(s)$ and $w$ be a weight of $\Fq$. By Sections \ref{section_computation_LLT_basis} and 
\ref{section_action_H}, we get a basis $\mathbb{B}$ of $\Fq \langle w \rangle$ that is 
$\barre$-invariant. Let $T(q)$ denote the transition matrix between the standard basis of $\Fq \langle w \rangle$ and the basis 
$\mathbb{B}$, and let $A(q)$ denote the matrix of the involution $\barre$ of $\Fq \langle w \rangle$ with
respect to the standard basis. Since the basis $\mathbb{B}$ is $\barre$-invariant, we have the following:
\begin{lemma} With the notation above, we have $A(q)=T(q)\big(T(q^{-1})\big)^{-1}$. \cqfd
\end{lemma}
Thus this lemma gives an algorithm for computing the involution $\barre$ of $\Fq \langle w \rangle$. Once this is known,
the computation of the canonical bases can be performed in a classical way
by solving unitriangular systems (see \eg \cite[Proof of Thm. 7.1]{L}).

\section{Examples} \label{section_tables} 

In this section, we give some transition matrices of the canonical bases of some weight subspaces of higher-level
Fock spaces. We first give, for $k \in \mZ$, the matrices

$$\Delta_k(q)=\bigl(\Delta^{+}_{\bd{\lambda}_l,\bd{\mu}_l ;\, \bd{s}_l^{(k)}}(q) \bigr)
_{\bd{\lambda}_l,\, \bd{\mu}_l \in \Pi^l(\bd{s}_l^{(k)} ;\, w^{(k)})}$$
for $n=l=2$, $\bd{s}_l^{(k)}=(2k,-2k)$ and $w^{(k)}=\wt(|\bd{\emptyset}_l,\bd{s}_l^{(k)} \rangle) - (2 \alpha_0+2 \alpha_1)$. These matrices should 
be read by columns, for example we have
$$G^{+}\Big( \big( (3,1),\emptyset \big), (-2,2) \Big) = \Big| \big( (3,1),\emptyset \big), (-2,2) \Big \rangle \,+\,
q\, \Big| \big( (2,2),\emptyset \big), (-2,2) \Big \rangle \,+\, q^2\,\Big| \big( (2,1,1),\emptyset \big), (-2,2) \Big \rangle.$$ 
Let $N_k$ denote the number of nonzero entries in $\Delta_k(q)$. Mark by a $*$ the rows indexed by the multi-partitions
in $\Pi^l(\bd{s}_l^{(k)},n)^{°}$. As a consequence, the corresponding columns give the expression (in the standard basis) of the vectors
of the lower global crystal basis of the irreducible $\Uq$-module $\Mq[\bd{s}_l^{(k)}]$. \\

Using Uglov's algorithm, we were able to compute the matrices $\Delta_k(q)$ only for $|k| \leq 1$, because 
otherwise the number of factors of the $q$-wedge products to be straightened was too large. With our algorithm
we computed the matrices $\Delta_k(q)$ for $|k| \lessapprox 20$. The matrix $\Delta_0(q)$ was already published in \cite{U}. Note that for a given $k \geq 1$, the matrices $\Delta_k(q)$ and $\Delta_{-k}(q)$ are equal up to a permutation of the rows and the columns; if $k \geq 2$, this is a special case of \cite[Prop. 5.5]{Y2}. Moreover, note that for all $k \geq 1$ (case of the dominant
multi-charges), the matrices $\Delta_k(q)$ are equal (up to a permutation of the rows and the columns). This property was observed and proved (for $k \geq 6$) in \cite[Thm 5.2]{Y2}. 
%!!!!!!!!!!!! references in [Y2]
This result supports in turn a conjecture for computing the $q$-decomposition matrices of some cyclotomic $v$-Schur algebras (see Section \ref{section_intro} and \cite{Y1}). In this conjecture, the assumption of dominance is necessary. Indeed, note that in our example we have $N_k \neq N_0$ for $k \neq 0$, 
which shows that $\Delta_k(1)$ $(k \neq 0)$ cannot be obtained from $\Delta_0(1)$
by a mere permutation of rows and columns. 

\begin{itemize}

\item[$*$] \underline{The matrices $\Delta_k(q)$ for $k \leq -1$}

$$\scriptsize \Delta_k(q)= \begin{array}{ll}
 \left( \begin{array}{*{16}{c}}
 1 & . & . & . & . & . & . & . & . & . & . & . & . & . & . & . \\
 q & 1 & . & . & . & . & . & . & . & . & . & . & . & . & . & . \\
 0 & q & 1 & . & . & . & . & . & . & . & . & . & . & . & . & . \\
 q & q^2 & q & 1 & . & . & . & . & . & . & . & . & . & . & . & . \\
 q^2 & q^3 & q^2 & q & 1 & . & . & . & . & . & . & . & . & . & . & . \\
 q & 0 & 0 & 0 & q & 1 & . & . & . & . & . & . & . & . & . & . \\
 q^2 & 0 & 0 & q & q^2 & q & 1 & . & . & . & . & . & . & . & . & . \\
 q^2 & 0 & 0 & q & 0 & 0 & 0 & 1 & . & . & . & . & . & . & . & . \\
 q^2 & 0 & 0 & q & q^2 & q & 0 & 0 & 1 & . & . & . & . & . & . & . \\
 q^3 & 0 & q & 2q^2 & q^3 & q^2 & q & q & q & 1 & . & . & . & . & . & . \\
 q^2 & q & q^2 & q^3 & q^4 & q^3+q & q^2 & 0 & q^2 & q & 1 & . & . & . & . & . \\
 q^2 & 0 & 0 & 0 & 0 & q & 0 & 0 & 0 & 0 & 0 & 1 & . & . & . & . \\
 q^3 & q^2 & 0 & 0 & 0 & q^2 & 0 & 0 & 0 & 0 & q & q & 1 & . & . & . \\
 0 & q^3 & q^2 & 0 & 0 & 0 & 0 & 0 & 0 & q & q^2 & 0 & q & 1 & . & . \\
 q^3 & q^4 & q^3 & q^2 & 0 & q^2 & q & 0 & q & q^2 & q^3 & q & q^2 & q & 1 & . \\
 q^4 & 0 & 0 & q^3 & 0 & q^3 & q^2 & q^2 & q^2 & q & 0 & q^2 & 0 & 0 & q & 1 \\
 \end{array} \right)
& \hspace{-5mm}
 \begin{array}{ll} 
 * & \hspace{-3mm} \bigl( \emptyset,(4) \bigr) \\ 
 * & \hspace{-3mm} \bigl( \emptyset,(3, 1) \bigr) \\
   & \hspace{-3mm} \bigl( \emptyset, (2, 2) \bigr) \\
   & \hspace{-3mm} \bigl( \emptyset, (2, 1, 1) \bigr) \\
 * & \hspace{-3mm} \bigl( (1), (2, 1) \bigr) \\ 
   & \hspace{-3mm} \bigl( (2), (2) \bigr) \\ 
   & \hspace{-3mm} \bigl( (1, 1), (2) \bigr) \\  
   & \hspace{-3mm} \bigl( \emptyset, (1, 1, 1, 1) \bigr) \\
   & \hspace{-3mm} \bigl( (2), (1, 1) \bigr) \\ 
   & \hspace{-3mm} \bigl( (1, 1), (1, 1) \bigr) \\ 
   & \hspace{-3mm} \bigl( (2, 1), (1) \bigr) \\ 
   & \hspace{-3mm} \bigl( (4),\emptyset \bigr) \\
   & \hspace{-3mm} \bigl( (3, 1),\emptyset \bigr) \\ 
   & \hspace{-3mm} \bigl( (2, 2), \emptyset \bigr) \\ 
   & \hspace{-3mm} \bigl( (2, 1, 1), \emptyset \bigr) \\ 
   & \hspace{-3mm} \bigl( (1, 1, 1, 1), \emptyset \bigr) \\ 
 \end{array}
\end{array},$$ 
$N_k=87$

\newpage %!!!!!!!!!!!!!!!!!
\item[$*$] \underline{The matrix $\Delta_0(q)$}

$$\scriptsize \Delta_0(q)= \begin{array}{ll}
 \left( \begin{array}{*{16}{c}}
1 & . & . & . & . & . & . & . & . & . & . & . & . & . & . & . \\
q & 1 & . & . & . & . & . & . & . & . & . & . & . & . & . & . \\
q & 0 & 1 & . & . & . & . & . & . & . & . & . & . & . & . & . \\
q^2 & q & q & 1 & . & . & . & . & . & . & . & . & . & . & . & . \\
0 & q & 0 & 0 & 1 & . & . & . & . & . & . & . & . & . & . & . \\
q^2 & 0 & q & 0 & 0 & 1 & . & . & . & . & . & . & . & . & . & . \\
q & q^2 & 0 & 0 & q & q & 1 & . & . & . & . & . & . & . & . & . \\
q^2 & 0 & q & 0 & 0 & q^2 & q & 1 & . & . & . & . & . & . & . & . \\
q^2 & q^3 & 0 & 0 & q^2 & 0 & q & 0 & 1 & . & . & . & . & . & . & . \\
q^2 & 0 & q & 0 & 0 & q^2 & q & 0 & 0 & 1 & . & . & . & . & . & .  \\
q^3 & q^2 & 2q^2 & q & q & q^3 & q^2 & q & 0 & q & 1 & . & . & . & . & . \\
q^2 & 0 & q^3 & 0 & q^2 & q^4 & q^3+q & q^2 & 0 & q^2 & q & 1 & . & . & . & . \\
q^3 & 0 & 0 & 0 & 0 & 0 & q^2 & 0 & q & 0 & 0 & q & 1 & . & . & . \\
0 & q^3 & 0 & q^2 & q^2 & 0 & 0 & 0 & 0 & 0 & q & 0 & 0 & 1 & . & . \\
q^3 & q^4 & q^2 & q^3 & q^3 & 0 & q^2 & q & q & q & q^2 & 0 & 0 & q & 1 & . \\
q^4 & 0 & q^3 & 0 & 0 & 0 & q^3 & q^2 & q^2 & q^2 & q & q^2 & q & 0 & q & 1  \\
 \end{array} \right)
& \hspace{-5mm}
 \begin{array}{ll} 
 * & \hspace{-3mm} \bigl( \emptyset,(4) \bigr) \\ 
 * & \hspace{-3mm} \bigl( \emptyset,(3, 1) \bigr) \\
   & \hspace{-3mm} \bigl( (4),\emptyset \bigr) \\
   & \hspace{-3mm} \bigl( (3, 1),\emptyset \bigr) \\ 
   & \hspace{-3mm} \bigl( \emptyset, (2, 2) \bigr) \\ 
 * & \hspace{-3mm} \bigl( (2), (2) \bigr) \\ 
   & \hspace{-3mm} \bigl( (1), (2, 1) \bigr) \\ 
   & \hspace{-3mm} \bigl( (1, 1), (2) \bigr) \\  
   & \hspace{-3mm} \bigl( \emptyset, (2, 1, 1) \bigr) \\ 
   & \hspace{-3mm} \bigl( (2), (1, 1) \bigr) \\ 
   & \hspace{-3mm} \bigl( (2, 1), (1) \bigr) \\ 
   & \hspace{-3mm} \bigl( (1, 1), (1, 1) \bigr) \\ 
   & \hspace{-3mm} \bigl( \emptyset, (1, 1, 1, 1) \bigr) \\ 
   & \hspace{-3mm} \bigl( (2, 2), \emptyset \bigr) \\ 
   & \hspace{-3mm} \bigl( (2, 1, 1), \emptyset \bigr) \\ 
   & \hspace{-3mm} \bigl( (1, 1, 1, 1), \emptyset \bigr) \\ 
 \end{array}
\end{array},$$
$N_0=86$

\item[$*$] \underline{The matrices $\Delta_k(q)$ for $k \geq 1$}

$$\scriptsize \Delta_k(q)= \begin{array}{ll}
 \left( \begin{array}{*{16}{c}}
 1 & . & . & . & . & . & . & . & . & . & . & . & . & . & . & . \\
 q & 1 & . & . & . & . & . & . & . & . & . & . & . & . & . & . \\
 0 & q & 1 & . & . & . & . & . & . & . & . & . & . & . & . & . \\
 q & q^2 & q & 1 & . & . & . & . & . & . & . & . & . & . & . & . \\
 q^2 & q^3 & q^2 & q & 1 & . & . & . & . & . & . & . & . & . & . & . \\
 q & 0 & 0 & 0 & q & 1 & . & . & . & . & . & . & . & . & . & . \\
 q^2 & 0 & 0 & q & q^2 & q & 1 & . & . & . & . & . & . & . & . & . \\
 q^2 & 0 & 0 & q & 0 & 0 & 0 & 1 & . & . & . & . & . & . & . & . \\
 q^2 & 0 & 0 & q & q^2 & q & 0 & 0 & 1 & . & . & . & . & . & . & . \\
 q^3 & 0 & q & 2q^2 & q^3 & q^2 & q & q & q & 1 & . & . & . & . & . & . \\
 q^2 & q & q^2 & q^3 & q^4 & q^3+q & q^2 & 0 & q^2 & q & 1 & . & . & . & . & . \\
 q^2 & 0 & 0 & 0 & 0 & q & 0 & 0 & 0 & 0 & 0 & 1 & . & . & . & . \\
 q^3 & q^2 & 0 & 0 & 0 & q^2 & 0 & 0 & 0 & 0 & q & q & 1 & . & . & . \\
 0 & q^3 & q^2 & 0 & 0 & 0 & 0 & 0 & 0 & q & q^2 & 0 & q & 1 & . & . \\
 q^3 & q^4 & q^3 & q^2 & 0 & q^2 & q & 0 & q & q^2 & q^3 & q & q^2 & q & 1 & . \\
 q^4 & 0 & 0 & q^3 & 0 & q^3 & q^2 & q^2 & q^2 & q & 0 & q^2 & 0 & 0 & q & 1 \\
 \end{array} \right)
& \hspace{-5mm}
 \begin{array}{ll} 
 * & \hspace{-3mm} \bigl( (4),\emptyset \bigr) \\
 * & \hspace{-3mm} \bigl( (3, 1),\emptyset \bigr) \\ 
   & \hspace{-3mm} \bigl( (2, 2), \emptyset \bigr) \\ 
   & \hspace{-3mm} \bigl( (2, 1, 1), \emptyset \bigr) \\ 
 * & \hspace{-3mm} \bigl( (2, 1), (1) \bigr) \\ 
   & \hspace{-3mm} \bigl( (2), (2) \bigr) \\
   & \hspace{-3mm} \bigl( (2), (1, 1) \bigr) \\ 
   & \hspace{-3mm} \bigl( (1, 1, 1, 1), \emptyset \bigr) \\ 
   & \hspace{-3mm} \bigl( (1, 1), (2) \bigr) \\
   & \hspace{-3mm} \bigl( (1, 1), (1, 1) \bigr) \\
   & \hspace{-3mm} \bigl( (1), (2, 1) \bigr) \\
   & \hspace{-3mm} \bigl( \emptyset,(4) \bigr) \\ 
   & \hspace{-3mm} \bigl( \emptyset,(3, 1) \bigr) \\
   & \hspace{-3mm} \bigl( \emptyset, (2, 2) \bigr) \\ 
   & \hspace{-3mm} \bigl( \emptyset, (2, 1, 1) \bigr) \\ 
   & \hspace{-3mm} \bigl( \emptyset, (1, 1, 1, 1) \bigr) \\ 
  \end{array}
\end{array},$$ 
$N_k=87$
\end{itemize}

%\bigskip

We now give the matrix $\Delta(q)=\bigl(\Delta^{+}_{\bd{\lambda}_l,\bd{\mu}_l ;\, \bd{s}_l}(q) \bigr)
_{\bd{\lambda}_l,\, \bd{\mu}_l \in \Pi^l(\bd{s}_l ;\, w)}$
for $n=l=2$, $\bd{s}_l=(0,0)$ and $w=\wt(|\bd{\emptyset}_l,\bd{s}_l \rangle) - (7 \alpha_0+4 \alpha_1)$.
It would be very hard to compute it using Uglov's algorithm. Indeed, the
partitions $\lambda$ such that $|\lambda,0 \rangle \in \Fq \langle w \rangle$ all satisfy $|\lambda|=25$, so 
applying Uglov's algorithm requires the straightening of $q$-wedge products with at least $25$ factors. Moreover,
since $\dim (\Fq \langle w \rangle) = 28$, the matrix $\Delta(q)$ here is too large to be displayed on a single page. We 
thus write it as 
$$\Delta(q)= \left( \begin{array}{cc} \Delta^{(1,1)}(q) & 0 \\ \Delta^{(2,1)}(q) & \Delta^{(2,2)}(q) \end{array} \right),$$
where $\Delta^{(1,1)}(q)$, $\Delta^{(2,1)}(q)$ and $\Delta^{(2,2)}(q)$ are $(14 \times 14)$-matrices. 

$$\scriptsize \Delta^{(1,1)}(q) = \begin{array}{ll}
 \left( \begin{array}{*{14}{c}}

1 & . & . & . & . & . & . & . & . & . & . & . & . & . \\
q & 1 & . & . & . & . & . & . & . & . & . & . & . & . \\ 
0 & 0 & 1 & . & . & . & . & . & . & . & . & . & . & . \\ 
q & 0 & q & 1 & . & . & . & . & . & . & . & . & . & . \\
q^2 & q & 0 & q & 1 & . & . & . & . & . & . & . & . & . \\ 
q & 0 & q & q^2 & q & 1 & . & . & . & . & . & . & . & . \\ 
q^2 & 0 & q^2 & 0 & 0 & q & 1 & . & . & . & . & . & . &. \\ 
q^2 & 0 & 0 & q & 0 & 0 & 0 & 1 & . & . & . & . & . & . \\ 
q^3 & q^2 & q & q^2 & q & 0 & 0 & q & 1 & . & . & . & . & . \\ 
0 & 0 & q^2 & 0 & 0 & 0 & 0 & 0 & q & 1 & . & . & . & . \\ 
q^2 & q & q^2 & q^3 & q^2 & q & 0 & q^2 & q & 0 & 1 & . & . & . \\ 
q^3 & q^2 & q^3 & 0 & q & q^2 & q & 0 & q^2 & q & q & 1 & . & . \\ 
0 & 0 & q^2 & q & 0 & 0 & 0 & 0 & 0 & 0 & 0 & 0 & 1 & . \\ 
0 & 0 & q^3 & q^2 & q & 0 & 0 & q & q^2 & q & 0 & 0 & q & 1 \\
 \end{array} \right)
& \hspace{-5mm}
 \begin{array}{ll} 
 * & \hspace{-3mm} \bigl( (1), (7, 2, 1) \bigr) \\
   & \hspace{-3mm} \bigl( (7, 2, 1), (1) \bigr) \\
 * & \hspace{-3mm} \bigl( (1), (5, 4, 1) \bigr) \\
 * & \hspace{-3mm} \bigl( (3), (5, 2, 1) \bigr) \\
   & \hspace{-3mm} \bigl( (3, 2, 1), (5) \bigr) \\
   & \hspace{-3mm} \bigl( (1, 1, 1), (5, 2, 1) \bigr) \\
   & \hspace{-3mm} \bigl( (1), (5, 2, 1, 1, 1) \bigr) \\
   & \hspace{-3mm} \bigl( (5), (3, 2, 1) \bigr) \\
   & \hspace{-3mm} \bigl( (5, 2, 1), (3) \bigr) \\
   & \hspace{-3mm} \bigl( (5, 4, 1), (1) \bigr) \\
   & \hspace{-3mm} \bigl( (5, 2, 1), (1, 1, 1) \bigr) \\
   & \hspace{-3mm} \bigl( (5, 2, 1, 1, 1), (1) \bigr) \\
 * & \hspace{-3mm} \bigl( (3, 2), (3, 2, 1) \bigr) \\
   & \hspace{-3mm} \bigl( (3, 2, 1), (3, 2) \bigr) \\
 \end{array}
\end{array}$$ \\

$$\scriptsize \hspace{-10mm} \Delta^{(2,1)}(q) = \hspace{-3mm} \begin{array}{ll}
 \left( \begin{array}{*{14}{c}}

q^2 & 0 & q^2 & q^3+q & q^2 & q & 0 & q^2 & 0 & 0 & 0 & 0 & q^2 & q \\
q^3 & 0 & q^3 & q^2 & 0 & q^2 & q & 0 & 0 & 0 & 0 & 0 & 0 & 0 \\
q^3 & q^2 & q^3 & q^4+q^2 & q^3+q & q^2 & 0 & q^3 & q^2 & q & q & 0 & q^3 & q^2 \\
q^4 & q^3 & q^4 & q^3 & 2q^2 & q^3 & q^2 & q^2 & q^3 & q^2 & q^2 & q & 0 & q \\
0 & 0 & q^2 & q^3 & q^2 & q & 0 & 0 & 0 & q^2 & 0 & 0 & q^4 & q^3 \\
0 & 0 & q^3 & 0 & 0 & q^2 & q & 0 & 0 & 0 & 0 & 0 & 0 & 0 \\
q^2 & 0 & q^4 & q^3 & q^2 & q^3+q & q^2 & 0 & 0 & 0 & 0 & 0 & 0 & 0 \\
q^3 & 0 & 0 & q^4 & q^3 & q^2 & 0 & q^3 & 0 & 0 & 0 & 0 & 0 & q^2 \\
q^4 & 0 & 0 & 0 & 0 & q^3 & q^2 & 0 & 0 & 0 & 0 & 0 & 0 & 0 \\
0 & 0 & q^3 & q^4 & q^3 & q^2 & 0 & q^3 & q^2 & q^3 & q & 0 & q^5 & q^4 \\
q^3 & q^2 & 0 & q^4 & q^3 & q^2 & 0 & 0 & 0 & 0 & q & 0 & 0 & 0 \\
q^4 & q^3 & q^4 & q^5 & q^4+q^2 & 2q^3 & q^2 & q^4 & q^3 & q^2 & 2q^2 & q & 0 & q^3 \\
0 & 0 & q^5 & 0 & q^3 & q^4 & q^3 & 0 & q^4 & q^3+q & q^3 & q^2 & 0 & 0 \\
q^5 & q^4 & 0 & 0 & q^3 & q^4 & q^3 & 0 & 0 & 0 & q^3 & q^2 & 0 & 0 \\ 
 \end{array} \right)
& \hspace{-5mm}
 \begin{array}{l}
 \bigl( (3, 1, 1), (3, 2, 1) \bigr) \\
 \bigl( (3), (3, 2, 1, 1, 1) \bigr) \\
 \bigl( (3, 2, 1), (3, 1, 1) \bigr) \\
 \bigl( (3, 2, 1, 1, 1), (3) \bigr) \\
 \bigl( (2, 2, 1), (3, 2, 1) \bigr) \\
 \bigl( (1), (3, 2, 2, 2, 1) \bigr) \\
 \bigl( (1, 1, 1), (3, 2, 1, 1, 1) \bigr) \\
 \bigl( (1, 1, 1, 1, 1), (3, 2, 1) \bigr) \\
 \bigl( (1), (3, 2, 1, 1, 1, 1, 1) \bigr) \\
 \bigl( (3, 2, 1), (2, 2, 1) \bigr) \\
 \bigl( (3, 2, 1), (1, 1, 1, 1, 1) \bigr) \\
 \bigl( (3, 2, 1, 1, 1), (1, 1, 1) \bigr) \\
 \bigl( (3, 2, 2, 2, 1), (1) \bigr) \\
 \bigl( (3, 2, 1, 1, 1, 1, 1), (1) \bigr) \\
 \end{array}
\end{array}$$ \\

$$\scriptsize \Delta^{(2,2)}(q) = \begin{array}{ll}
 \left( \begin{array}{*{14}{c}}
1 & . & . & . & . & . & . & . & . & . & . & . & . & . \\
q & 1 & . & . & . & . & . & . & . & . & . & . & . & . \\
q & 0 & 1 & . & . & . & . & . & . & . & . & . & . & . \\
q^2 & q & q & 1 & . & . & . & . & . & . & . & . & . & . \\
q^2 & 0 & q & 0 & 1 & . & . & . & . & . & . & . & . & . \\
0 & 0 & 0 & 0 & q & 1 & . & . & . & . & . & . & . & . \\
q^2 & q & q & 0 & q^2 & q & 1 & . & . & . & . & . & . & . \\
q^3 & q^2 & q^2 & q & 0 & 0 & q & 1 & . & . & . & . & . & . \\
0 & q & 0 & 0 & 0 & q & q^2 & q & 1 & . & . & . & . & . \\
q^3 & 0 & q^2 & 0 & q & 0 & 0 & 0 & 0 & 1 & . & . & . & . \\
q^3 & q^2 & q^2 & 0 & 0 & 0 & q & 0 & 0 & 0 & 1 & . & . & . \\
q^4 & q^3 & q^3+q & q^2 & q^2 & q & q^2 & q & 0 & q & q & 1 & . & . \\
0 & 0 & q^2 & 0 & q^3 & q^2 & 0 & 0 & 0 & q^2 & 0 & q & 1 & . \\
0 & q^2 & q^2 & q & 0 & q^2 & q^3 & q^2 & q & 0 & q^2 & q & 0 & 1 \\ 
 \end{array} \right)
& \hspace{-5mm}
 \begin{array}{l}
 \bigl( (3, 1, 1), (3, 2, 1) \bigr) \\
 \bigl( (3), (3, 2, 1, 1, 1) \bigr) \\
 \bigl( (3, 2, 1), (3, 1, 1) \bigr) \\
 \bigl( (3, 2, 1, 1, 1), (3) \bigr) \\
 \bigl( (2, 2, 1), (3, 2, 1) \bigr) \\
 \bigl( (1), (3, 2, 2, 2, 1) \bigr) \\
 \bigl( (1, 1, 1), (3, 2, 1, 1, 1) \bigr) \\
 \bigl( (1, 1, 1, 1, 1), (3, 2, 1) \bigr) \\
 \bigl( (1), (3, 2, 1, 1, 1, 1, 1) \bigr) \\
 \bigl( (3, 2, 1), (2, 2, 1) \bigr) \\
 \bigl( (3, 2, 1), (1, 1, 1, 1, 1) \bigr) \\
 \bigl( (3, 2, 1, 1, 1), (1, 1, 1) \bigr) \\
 \bigl( (3, 2, 2, 2, 1), (1) \bigr) \\
 \bigl( (3, 2, 1, 1, 1, 1, 1), (1) \bigr) \\
 \end{array}
\end{array}$$ \\

\vspace{5mm}
\small \rm Xavier YVONNE, Institut Camille Jordan (Math\'ematiques), Universit\'e Lyon I, 43 Bd du 11~novembre 1918, 69622 Villeurbanne Cedex, France. \\

\indent \it E-mail address: \tt yvonne@math.univ-lyon1.fr
\end{document}